\documentclass[12pt]{article}
\usepackage[dvipsnames,usenames]{color}
\usepackage{amsmath}
\usepackage{amsfonts}
\usepackage{amssymb}
\usepackage{mathrsfs}
\usepackage{amsthm}
\usepackage{mathrsfs}
\usepackage{amsmath}
\usepackage{amsfonts}
\usepackage{amssymb}
\usepackage{enumerate}

\usepackage{amsmath, amsthm, amsfonts, amssymb}
\def\<{\langle}
\def\>{\rangle}

\def\p{\partial}
\def\a{\alpha}

\def\O{\Omega}

\def\-{\overline}

\def\endpf{\hbox{\vrule height1.5ex width.5em}}

\def\b{\beta}
\def\a{\alpha}
\def\endpf{\hbox{\vrule height1.5ex width.5em}}

\def\-{\overline}

\def\O{\Omega}

\def\p{\partial}

\def\endpf{\hbox{\vrule height1.5ex width.5em}}

\def\a{\alpha}

\def\Ol{\overline}

\def\b{\beta}

\def\a{\alpha}

\def\endpf{\hbox{\vrule height1.5ex width.5em}}
\def\a{\alpha}

\def\ssubset{\subset\joinrel\subset}

\def\p{\partial}

\def\-{\overline}

\newtheorem{theorem}{Theorem}[section]
\newtheorem{lemma}[theorem]{Lemma}
\newtheorem{corollary}[theorem]{Corollary}
\newtheorem{proposition}[theorem]{Proposition}

\newtheorem{definition}[theorem]{Definition}
\newtheorem{example}[theorem]{Example}
\newtheorem*{examplenonum}{Example}
\newtheorem{remark}[theorem]{Remark}

\date{}

\begin{document}

\title{\bf Bergman-Einstein metrics, hyperbolic metrics and Stein spaces with spherical
boundaries}

\author{Xiaojun Huang \footnote{Supported in part by NSF grant  DMS-1665412.}\  \ and  Ming Xiao
\footnote{Supported in part by NSF grant DMS-1800549}}

\vspace{3cm} \maketitle

\begin{abstract} We give an affirmative solution to a conjecture of Cheng proposed in 1979
which asserts that the Bergman metric of a smoothly bounded strongly
pseudoconvex domain in $\mathbb{C}^n, n\geq 2,$ is K\"ahler-Einstein
if and only if the domain is biholomorphic to the ball. We establish
versions of various classical theorems that are used in the solution
for Stein spaces. Among other things, we construct  a hyperbolic
metric over a Stein space with spherical boundary. We also prove the
Q. K. Lu type uniformization  theorem for Stein spaces with isolated
normal singularities.

\medskip


\end{abstract}

\section{Introduction}
Canonical metrics are important objects under study in Complex
Analysis of Several Variables. Since Cheng-Yau  proved in \cite{CY} the
existence of a complete K\"ahler-Einstein metric over a bounded
pseudoconvex domain in ${\mathbb C}^n$ with reasonably smooth
boundary, it has become a natural question to understand when the
Cheng-Yau metric of a bounded pseudoconvex domain  is precisely its
 Bergman metric.
S. Y.  Cheng conjectured in 1979 \cite{Ch} that if the Bergman
metric of a smoothly bounded strictly pseudoconvex domain is
K\"ahler-Einstein, then the domain is biholomorphic to the ball.
Cheng's conjecture was previously obtained  by Fu-Wong  \cite{FW} and
Nemirovski-Shafikov \cite{NS1} in the case of complex dimension two.
 There are also  closely related studies on
versions of the Cheng conjecture in terms of metrics defined by
other canonical potential functions. The reader is referred to  work
of Li \cite{L1},\cite{L2}, \cite{L3} and the references therein on
this matter

The main purpose of this paper is twofold. One is to present a
solution of the Cheng conjecture in any dimensions. The other is to
use this opportunity to generalize several  classical results used
for the solution of the Cheng conjecture to Stein spaces with
singularities, which might be of independent interest and importance
in their own right.

The paper is divided into three parts.   In the first part of the
paper, we give an affirmative solution to Cheng's conjecture
\cite{Ch}, based on deep works of many mathematicians in the past 40
years. Namely, we prove the following:


\begin{theorem}\label{theorem1}
The Bergman metric of a smoothly bounded strongly pseudoconvex
domain in ${\mathbb C}^n$  $(n\ge 2)$ is K\"ahler-Einstein if and
only if the domain is biholomorphic to the ball.
\end{theorem}

To verify the Cheng conjecture, we first show that the Einstein
property of the Bergman metric over a bounded strongly pseudoconvex
domain $\Omega$ forces the boundary $\partial\Omega$ to be
spherical. Namely, at each point of $\partial\Omega$ there is a
small open piece of $\partial\Omega$ that is CR-diffeomorphic to an
open piece of the sphere of the same dimension. To prove that, we
will fundamentally  make use of the work done in Chern-Moser
\cite{CM}, Fefferman \cite{Fe1}\cite{Fe2}\cite{Fe3}, Christoffers
\cite{C}, Fu-Wong \cite{FW}, etc. Once this is known, as in the work
of  Nemirovski-Shafikov \cite{NS2}, one  can use the classical
Kerner theorem \cite{K} to prove that $\Omega$ is a ball quotient
(This can also be proved by applying  the extension theorem in
Chern-Ji \cite{CJ} instead of the general Kerner theorem, once
$\Omega$ is known to have a spherical boundary). In the second part
of the paper, we establish the following Kerner type theorem for
Stein spaces even with isolated complex singularities. Indeed, we
will prove a more general version of the multi-valued Hartogs
extension result, (see Theorem \ref{thm33}), which gives Theorem
\ref{thm33-preliminary} as a special case. See Section 3 for the
explanation of terminologies used below.

\begin{theorem}\label{thm33-preliminary}
Let $\Omega$ be a Stein space (of complex dimension at least two)
with possibly isolated singularities and connected compact strongly
pseudoconvex smooth boundary $M=\partial \Omega.$ Let $(f, D)$ be a non-constant
CR mapping, where $D$ is an open connected subset of $M$.  Let $C> 0$ be a constant such that $\|f\|\le C$ over
$D$.
Suppose that $(f, D)$ admits a  CR continuation along any curve in
$M$ and for each CR mapping element $(g, D^*)$ with $D^*\subset M$
obtained by continuing $(f,D)$ along a curve in $M$, it holds that
$\|g\|\le C$. Then $(f, D)$ admits a
holomorphic continuation along any curve $\gamma$ with $\gamma
(0)\in D$ and $\gamma (t)\in \hbox{Reg}(\Omega)$ for $t\in (0,1]$.
Moreover, for any holomorphic mapping element $(h, U)$ with
$U\subset \hbox{Reg}(\Omega)$ obtained from continuation of $(f,D)$,
we have $\|h\|<C$ on $U$.
\end{theorem}

Here $\mathrm{Reg}(\Omega)$ denotes the set of smooth points in
$\Omega.$ We mention that in the above theorem we do not have in
general the extension up to the singular points even if the
singularities are normal, which is different from  the classical
Hartogs extension theorem. (See Remark \ref{rekk}). Kerner's
original theorem and his method of proof \cite {K} do not apply to
Stein spaces with singularities. Also in the work of Chern-Ji
\cite{CJ} on the spherical boundary case with $(f,D)$ being a CR
differomorphism to a piece of the sphere, the lift of $(f,D)$ was
proved to extend along any curve in the extended G-structure bundle
over $\Omega$ which implies the bimeromorphic extension of $(f,D)$
on the base manifold $\Omega$ along any cure. In their proof, it is
crucial for $\Omega$ to be a domain in ${\mathbb C}^n$. For the
proof of Theorem \ref{thm33} and thus Theorem
\ref{thm33-preliminary}, we will employ a different but in fact more
elementary and self-contained argument than those used in \cite{K}
and \cite{CJ}, based on the Lewy and Baouendi-Treves extension
theorem \cite{BER}, Morse function theory and the
Phragm\'en-Lindel\"off principle. Making use of Theorem
\ref{thm33-preliminary}, we can construct a hyperbolic metric over
$\Omega$:


\begin{theorem}\label{theorem2}
Let $\Ol{\Omega}$ be a Stein space of dimension at least two that
has  a connected compact  smooth boundary and possibly isolated singularities
in the interior. Assume the boundary  $\partial\Omega$ is spherical.
Then there is a hyperbolic metric $\omega_0$ on $\mathrm{Reg}(\Omega)$ that is complete at infinity.
\end{theorem}

We mention that Theorem \ref{theorem2}  is only new when $\Omega$
has singularities. Also, a K\"ahler metric is called hyperbolic if
it has constant negative holomorphic sectional curvature. We say the
metric is complete at infinity if for every real number $R> 0,$ the
ball centered at $p_0 \in \mathrm{Reg}(\Omega)$ with radius $R$
(with respect to the metric) has compact closure in $\Omega.$

When $\Omega$ in Theorem \ref{theorem2} is not smooth, the
hyperbolic metric may not be complete at the singular point. For
instance, as proved in Huang \cite{Hu}, if $\Ol{\Omega}$ is embedded
in a complex Euclidean space with $\partial \Omega$ spherical and
algebraic, then ${\Omega}$ has exactly one singular point which is a
finite quotient singularity of the unit ball $\mathbb{B}^n \subset
\mathbb{C}^n$. The naturally inherited hyperbolic metric is then not
complete at the singular point. On the other hand, if the hyperbolic
metric in Theorem \ref{theorem2} is complete at the singular point,
then $\Omega$ must be a completion of a non-compact ball quotient
space. A typical example of such a domain is obtained by taking the
quotient space of $\mathbb{B}^n$ for a lattice $\Gamma \subset
Aut({\mathbb B}^n)$ with a fundamental region having  one boundary
component shrinking to a point on the unit sphere and one boundary
component being an open piece of the unit sphere. The possible
non-completeness at the singularities makes it subtle to study the uniqueness of such
hyperbolic metrics. The following Example \ref{epdisk} shows that
the hyperbolic metric in Theorem \ref{theorem2} may not be unique
(even up to scaling) in general in the dimension one case. When the
dimension is at least two, we will establish a  uniqueness result of
the metric under a geodesic completeness assumption near the
boundary. (See Proposition \ref{propcomplete}).

\begin{example}\label{epdisk}
Let $\Delta$ be the unit disk and let $X$ be the Riemann surface
with singularity in $\mathbb{C}^2$ give by $\{ (z, w) \in
\mathbb{C}^2: w^2=z^3\} \cap \Delta^2.$ It is indeed the image of
the map $t \rightarrow (t^2, t^3), t \in \Delta$, and has an
isolated singularity $(0, 0).$ Note $X^*:=X \setminus \{ (0, 0) \}$
is biholomorphic to the punctured disk $\Delta^*:=\Delta \setminus
\{ 0 \}.$ Then the hyperbolic metric on $\Delta^*$ induces a
hyperbolic metric $\omega_1$ on $X^*$ and the Bergman metric on
$\Delta^*$ induces a  hyperbolic metric $\omega_2$ on $X^*$. Notice
that $\omega_1$ is complete near the singularity while $\omega_2$ is
not. And they are both complete at infinity. A scaling cannot make
them identical.
\end{example}


Back to our description of the proof of Theorem \ref{theorem1}, once
the domain $\Omega$  is known to be a ball quotient (since a domain
in ${\mathbb C}^n$ has no singularities), one can apply a classical
theorem of Qi-Keng Lu to conclude that $\Omega$ is biholomorphic to
the unit ball. The last part of our paper is devoted to establishing
 the Lu type theorem for a Stein space with possibly normal singularities:

\begin{theorem}\label{theorem3}
Let $\Omega$ be a Stein space with possibly isolated normal
singularities of complex dimension $n \geq 2$ with Bergman metric
$\omega^B.$ Assume that the Bergman metric is complete at infinity
and the induced Bergman-Bochner map is one-to-one (over
$Reg(\Omega)$). Then $\Omega$ is biholomorphic to the unit ball
$\mathbb{B}^n$ if and only if $(\mathrm{Reg}(\Omega), \omega^B)$ has
constant holomorphic sectional curvature.
\end{theorem}

Here the Bergman-Bochner map is a map from the smooth part of
$\Omega$ to ${\mathbb P}^\infty$ through an orthonormal basis of the
Bergman space of $\Omega$. See $\S 4$ for more details of its
definition.


\section{Proof of Theorem \ref{theorem1}}

Let $\Omega=\{ z \in \mathbb{C}^n: \rho(z) > 0 \}$ be a strictly
pseudoconvex domain with a smooth defining function $\rho.$ In
\cite{Fe1}, Fefferman showed that the Bergman kernel function $K(z)=K(z,
\overline{z})$ of $\O$  has the asymptotic expansion
$$K(z)=\frac{\phi(z)}{\rho^{n+1}(z)}+ \psi(z)\log \rho(z),$$
where $\phi, \psi \in C^{\infty}(\overline{\Omega})$ and
$\phi|_{\p\Omega} \neq 0.$ In particular, if the boundary $\p\Omega$
of $\Omega$ is spherical, then $\psi$ vanishes to infinite order at
the boundary $\p\Omega.$

We first recall the notion of Fefferman defining functions or
Fefferman approximate solutions.
Consider the following Monge-Amp\`ere type equation introduced  in
\cite{Fe2}:
$$J(u):=(-1)^n \hbox{ det}\left(
                             \begin{array}{cc}
                               u & u_{\overline{\beta}} \\
                               u_{\alpha} & u_{\alpha\overline{\beta}} \\
                             \end{array}
                           \right) =u^{n+1}\hbox{det}((\log\frac{1}{u})_{\a\-\b})=1 ~\text{in}~\Omega,
$$
with $u=0$ on $b \Omega.$ Fefferman proved that for any bounded
strictly pseudoconvex domain $\Omega$ with smooth boundary,  there
is a smooth defining function $r$ of $\Omega$ such that
$J(r)=1+O(r^{n+1}),$ which  is called a Fefferman approximate
solution or a Fefferman defining function of $\O$. Moreover, if
$r_1,\ r_2$ are two Fefferman approximate solutions, then
$r_1-r_2=O(\rho^{n+2})$, where $\rho$ is a given defining function
of $\Omega$.

We next recall the  Moser normal form theory \cite{CM} and the notion of
Fefferman scalar boundary invariants (cf. \cite{Fe3}, \cite{G}): Let $M
\subset \mathbb{C}^n$ be a real analytic strictly pseudoconvex
hypersurface containing $p \in \mathbb{C}^n$ (Assume here
that $n\ge 3$). There exists a coordinates system
$(z,w):=(z_{1},...,z_{n-1},w)$ such that in the new coordinates,
$p=0$ and $M$ is defined  near $p$ by an equation of the form:
\begin{equation}\label{normM}
u=|z|^2 + \sum_{|\alpha|, |\beta| \geq 2, l \geq 0} A_{\alpha \overline{\beta}}^l z^{\alpha} \overline{z}^{\beta} v^l,
\end{equation}
where $w=u+iv, \alpha=(\alpha_{1},...,\alpha_{n}),
\beta=(\beta_{1},...,\beta_{n})$ are multiindices.  Moreover, the
coefficients $A_{\alpha \overline{\beta}}^l  \in \mathbb{C}$ satisfy
the following:

\begin{itemize}
\item $A_{\alpha \overline{\beta}}^l$ is symmetric with respect to permutation of the indices in $\alpha$ and $\beta,$ respectively.

\item $\overline{A_{\alpha \overline{\beta}}^l}=A_{\beta \overline{\alpha}}^l.$

\item $\mathrm{tr} A_{2\overline{2}}^l =0, \mathrm{tr}^2 A_{2\overline{3}}^l=0, \mathrm{tr}^3 A_{3\overline{3}}^l =0,$ where $A_{p\overline{q}}^l$ is the symmetric
    tensor $[A_{\alpha \overline{\beta}}^l]_{|\alpha|=p, |\beta|=q}$ on $\mathbb{C}^{n-1}$ and the
    traces are the usual tensorial traces with respect to $\delta_{i\overline{j}}.$

\end{itemize}

Here (\ref{normM}) is called a normal form of $M$ at $p.$ When $M$
is merely smooth, the expansion is in the formal sense.
 $[A_{\alpha \overline{\beta}}^l]$ are called the normal form coefficients.
 Recall that a boundary scalar invariant at $p\leftrightarrow 0$, or briefly
an invariant of weight $w \geq 0$, is a polynomial $P$ in the normal
form coefficients $[A_{\alpha \overline{\beta}}^l]$ of $\p \Omega$
satisfying certain transformation laws. (See \cite{Fe3} and \cite{G} for more
details on this transformation law). 
Using a Fefferman defining function in the asymptotic expansion of
the Bergman kernel function:
\begin{equation}\label{eqKr}
K(z)=\frac{\phi(z)}{r^{n+1}(z)}+ \psi(z)\log r(z),
\end{equation}
with $\phi, \psi \in C^{\infty}(\overline{\Omega}), \
\phi|_{\p\Omega} \neq 0,$ then $\phi~ \mathrm{mod}~r^{n+1},
\psi~\mathrm{mod}~r^{\infty}$ are locally determined.
 Moreover, if $\p \Omega$ is in its normal form at $p=0 \in b \Omega,$
 then any Taylor coefficient at $0$ of $\phi$ of order $\leq n,$ and that  of $\psi$ of any
 order
 is a universal polynomial in the normal form coefficients $[A_{\alpha \overline{\beta}}^l].$
 (See Boutet-Sj\"ostrand \cite{BS} and a related argument in \cite{Fe3}.) In particular,
  we state the following result from \cite{C}. (See also \cite{G}):
\begin{proposition}\label{lemmaC}(\cite{C}, \cite{G})
  Let $\Omega$ be as above and suppose that $\p \O$ is in the Moser
normal form up to sufficiently high order. Let $r$ be a Fefferman
defining function, and let $\phi,\ \psi$ be as in (\ref{eqKr}). Then
$\phi|_{\p \O}=\frac{n!}{\pi^n}$,\ $\phi=\frac{n!}{\pi^n}+O(r^2)$
and $P_{2}= \frac{\phi-\frac{n!}{\pi^n}}{r^2}|_{\p \Omega}$ defines
an invariant of weight $2$ at $0$. Furthermore, if $n=2,$ then
$P_{2}=0.$ If $n \geq 3, P_{2}=c_{n}||A_{2\overline{2}}^{0}||^2$ for
some universal constant $c_{n} \neq 0.$
\end{proposition}

 As mentioned earlier,
Theorem \ref{theorem1} is known
 in the case of $n=2$ in \cite{FW} and \cite{NS1}. We
next assume that $n \geq 3$.

\bigskip

{\it Proof of Theorem \ref{theorem1}:} Recall the Fefferman
asymptotic expansion:
\begin{equation}\label{expansion}
K(z)=\frac{\phi(z)}{\rho^{n+1}(z)}+\psi(z) \log \rho(z)=
\frac{\phi+\rho^{n+1}\psi
 \log \rho}{\rho^{n+1}}\ \ \ \hbox{for}\  z \in \Omega
\end{equation}
with $\phi, \psi \in C^{\infty}(\overline{\Omega})$ and $\phi|_{\p
\Omega} \neq 0,$ where $\rho \in C^{\infty}(\overline{\Omega})$ is a
smooth defining function of $\Omega$ with $\Omega=\{ z \in
\mathbb{C}^n: \rho(z) >0 \}.$ Since  $K(z) >0$ for  $z \in \Omega,$
we have
\begin{equation}\label{eqnomega}
\phi+ \rho^{n+1}\psi \log \rho > 0\ ~\text{for }~ z \in \Omega.
\end{equation}
Thus
\begin{equation}\label{eqnr0}
(K)^{-\frac{1}{n+1}}(z)=\frac{\rho}{(\phi+\rho^{n+1} \psi \log
\rho)^{\frac{1}{n+1}}}
\end{equation}
is well-defined in $\Omega.$

Let $\O$ be a smoothly bounded strongly pseudoconvex domain. We
notice that the K\"ahler-Einstein condition of the Bergman metric is
equivalent to the fact that $\log{K(z)}$ is a K\"ahler-Einstein
potential function of $\O$.   More precisely, we have
$J\left[\left(\frac{\pi^n}{n!}K(z)\right)^{-\frac{1}{n+1}}\right]=1$
for  $z \in \Omega.$ (See \cite{FW}).
Let $r_{0}(z):=\left(\frac{\pi^n}{n!}K\right)^{-\frac{1}{n+1}}.$
We hence have that  $r_{0}(z) >0$ and   $J(r_{0}) =1$ in $\Omega.$
  We
next recall the following result of Fu-Wong \cite{FW}:

\begin{proposition}
Let $\Omega=\{ z \in \mathbb{C}^n: \rho >0\}$ be a bounded strongly
pseudoconvex domain with a smooth defining function $\rho.$ If the
Bergman metric of $\Omega$ is K\"ahler-Einstein, then the
coefficient of the logarithmic term in Fefferman's expansion
(\ref{expansion}) vanishes to infinite order at $b \Omega,$ i.e.,
$\psi=O(\rho^k)$ for any $k > 0.$
\end{proposition}

As a consequence, $\phi+ \rho^{n+1}\psi \log \rho$ extends smoothly
to a neighborhood of $\overline{\Omega}.$ Since $\phi|_{\p \Omega}
\neq 0,$ we have
$$\phi+ \rho^{n+1}\psi \log \rho > 0~\text{for all}~z \in \overline{\Omega}.$$
Hence  $r_{0}$ extends smoothly to a neighborhood of
$\overline{\Omega}$  and it is then easy to conclude that $r_{0}$ is
a Fefferman defining function of $\Omega.$  Then from the way $r_0$
was constructed,  it follows that
\begin{equation} \label{001}
K(z)=\frac{n!}{\pi^n} r_{0}^{-(n+1)}.
\end{equation}

Comparing  (\ref{001})  with  (\ref{eqKr}), we arrive at the
conclusion that if we let $r=r_{0}$
 in (\ref{eqKr}), then $\phi \equiv \frac{n!}{\pi^n}.$ Then it follows from Proposition
  \ref{lemmaC} that $P_{2}=c_{n}||A_{2\overline{2}}^0||^2=0$ at $p\in \p \O$ if $\p \O$ is in the
   Moser normal form up to sufficiently
  high order at $p$ with $A_{2\overline{2}}^0$ being the Chern-Moser-Weyl-tensor at $p$. Consequently,
  $A_{2\overline{2}}^0=0$ in each Moser normal coordinates at each $\p \Omega$, for $c_{n} \neq 0.$ That is, every boundary point of $\p \O$
   is a CR umbilical point.
By a classical result of Chern-Moser,  $\p \Omega$ is spherical.  We
then apply Theorem \ref{theorem2} to conclude that $\Omega$ has a
complete hyperbolic metric $\omega_0$.  As in \cite{NS1},  by the uniqueness of
result of Cheng and Yau \cite{CY},  the Bergman metric of $D$ is
proportional to the metric $\omega_0$ and hence has constant
holomorphic sectional curvature. Finally it follows from a well-know
result due to Qi-Keng Lu \cite{Lu} or Theorem \ref{theorem3}  that $\O$ is
biholomorphic to the unit ball. $\endpf$

\section{Proof of Theorem \ref{thm33-preliminary} and \ref{theorem2}}

In this section, we give a proof for Theorem \ref{thm33-preliminary} and \ref{theorem2}. Recall
that $\Omega$ is called a Stein space with possibly isolated
singularities if $\Omega$ is a Stein space and its singular set
$\mathrm{sing}(\Omega)$ consists of only isolated points if it is
not empty. Let $\Omega$ be such a Stein space with a compact smooth
boundary $M.$  In the following context, we write
$\Ol{\Omega}=\Omega \cup M,$ $\mathrm{Reg}(\Omega)=\Omega \setminus
\mathrm{sing}(\Omega)$ and $\mathrm{Reg}(\Ol{\Omega})=\Ol{\Omega}
\setminus \mathrm{sing}(\Omega)$.

We next introduce some standard terminalogy in the study of extensions of CR and holomorphic functions.

We say $(f, D)$ is a (continuous) CR function element on $M$ if $D$
is a  connected open piece of $M$ and $f: D \rightarrow
\mathbb{C}^N$ is a continuous CR map. Similarly, we say $(g, U)$ is
a holomorphic function element on $\mathrm{Reg}(\overline{\Omega})$ if $U$ is a  connected open subset of
$\mathrm{Reg}(\overline{\Omega})$ and $g: U
\rightarrow \mathbb{C}^N$ is a continuous map that is holomorphic in
$U \cap \Omega.$ Note in this article, a "function" can be vector-valued. We will also sometimes use the notion of germs of functions. Fix a CR function element $(f, D)$ on $M$ and $p \in D$. We write
$(f)_p$ for the germ of CR function $f$ at $p$ on $M$. Similarly, for a holomorphic function element $(g, U)$ and $q \in U,$ we denote by $[g]_q$ the germ of holomorphic function $g$ at $q$.

Fix a CR function element $(f, D)$ on $M$ and a curve $\sigma: [0, 1] \rightarrow \mathrm{Reg}(\overline{\Omega})$ such that $\sigma(0) \in D$ (resp. a curve $\sigma: [0, 1] \rightarrow M$  such that $\sigma(0) \in D.$)  We make the following definition.

\begin{definition}
We say $(f, D)$ admits a holomorphic continuation along $\sigma$ (resp. admits a (continuous) CR continuation along $\sigma$) if there exists a collection of holomorphic function elements $\{ (f_j, U_j)\}_{j=0}^{k}$ on $\mathrm{Reg}(\overline{\Omega})$ (resp. a collection of CR function elements $\{(f_j, U_j)\}_{j=0}^{k}$ on M) such that
\begin{enumerate}[(1)]
\item $f_0=f$ near $\sigma(0)$ in $U_0 \cap D,$
\item There is a partition $0 = t_0 < t_1 < \cdots < t_{k+1} =1$ such that $\sigma ([t_j, t_{j+1}]) \subset U_j,$ for all $0 \leq j \leq k;$

\item $f_j =f_{j+1}$ on $U_j \cap U_{j+1}$ for $0 \leq j \leq k-1.$
\end{enumerate}

Here $(f_k, U_k)$ is called a branch of $(f, D)$ obtained by holomorphic (resp. CR) continuation along $\sigma$.
Fix a plurisubharmonic function $\psi: \mathbb{C}^N \rightarrow \mathbb{R}$ such that $\psi(f) \leq 0$ on $D.$ We further say $(f, D)$ admits a holomorphic continuation with $\psi-$estimate along $\sigma$ (resp. admits a CR continuation with $\psi-$estimate along $\sigma$ ) if there is a choice of $\{(f_j, U_j)\}_{j=0}^k$ such that $\psi(f_j) \leq 0$ on $U_j$ for all $j$.
\end{definition}

\begin{definition}\label{defn32}
Let $(f, D)$ and $\psi$ be as above, and $\hat{\Omega}$ an open connected subset of $\mathrm{Reg}(\overline{\Omega})$
containing $D$. We say $(f, D)$ admits unrestricted holomorphic continuation
in $\hat{\Omega}$ (resp. admits unrestricted CR continuation in $M$) with $\psi-$estimiate if $(f, D)$
 admits a holomorphic continuation with $\psi-$estimate along every curve $\sigma$ in $\hat{\Omega}$ with  $\sigma(0) \in D$
 (resp. if $(f, D)$ admits a CR continuation with $\psi-$estimate along every curve $\gamma$ with $\gamma(0) \in D$ in $M$).
\end{definition}

 In particular, we say $(f, D)$ admits unrestricted bounded holomorphic continuation in $\hat{\Omega}$ if in Definition \ref{defn32} $\psi$ can be taken to be $\psi(z)=||z||^2-C$ for some $C>0.$ We define unrestricted bounded CR continuation in a similar way.

To make notations simple, we will say the function element $(g, U)$ is a holomorphic branch of $(f, D)$ in $\hat{\Omega}$ if (a). $U \subset \hat{\Omega},$ and (b). There exists a path $\sigma$ in $\hat{\Omega}$ with $\sigma(0) \in D$ and $\sigma(1) \in U$ such that $(g, U)$ can be obtained by holomorphic continuation of $(f, D)$ along $\sigma.$ Similarly, we say $[g]_q$ is a holomorphic branch of $(f, D)$ in $\hat{\Omega}$ if there exists a holomorphic branch $(g, U)$ of $(f, D)$ in $\hat{\Omega}$ in the above sense and $q \in U.$
\smallskip

To prove Theorem \ref{theorem2}, we will need the following  result
which establishes the Kerner type theorem \cite{K} for the Stein
space with possible singularities.  It also contains Theorem
\ref{thm33-preliminary} as its special case.
\begin{theorem}\label{thm33}
Let $\Omega$ be a Stein space with possibly isolated singularities
and connected compact strongly pseudoconvex boundary  $M=\partial
\Omega.$ Let $(f, D)$ be a CR function (or map)  element on $M$ and
$\psi: \mathbb{C}^N \rightarrow \mathbb{R}$ a plurisubharmonic
function such that  $\psi(f) \leq 0$ on $D.$ Then the following
conclusions hold.
\begin{enumerate}[(1)]
\item Suppose that $\mathrm{dim}_{\mathbb{C}}(\Omega) \geq 3$ and $(f, D)$ admits unrestricted CR continuation in $M$ with $\psi-$estimate.  Then $(f, D)$ admits unrestricted holomorphic continuation in $\mathrm{Reg}(\overline{\Omega})$ with $\psi-$estimate.

\item Suppose that $\mathrm{dim}_{\mathbb{C}}(\Omega) =2$ and $(f, D)$ admits unrestricted bounded CR continuation in $M$ with $\psi-$estimate.  Then $(f, D)$ admits unrestricted bounded holomorphic continuation in $\mathrm{Reg}(\overline{\Omega})$ with $\psi-$estimate.
\end{enumerate}

Moreover, in both cases, assume there is a holomorphic branch $(h,
U)$ of $(f, D)$ in $\mathrm{Reg}(\Ol{\Omega})$  such that $\psi
(h(p))=0$ at some point $p$ in $U \setminus M$.
Then $\psi(g) \equiv 0$ for any holomorphic branch $(g, V)$ of $(f, D).$ In particular, $\psi(f) \equiv 0$ on $D$.

\end{theorem}

The last part of Theorem \ref{thm33} immediately implies that if we additionally assume $\psi(f)$ is not constantly zero on $D$, then for any branch $(g, V)$ of $(f, D)$ in $\mathrm{Reg}(\Omega),$ it holds that $\psi(g) < 0$ on $V.$ Besides this, we also have the following consequence.

\begin{remark}\label{remarkstrict}
Assume in Theorem \ref{thm33} that $\psi$ is strictly plurisubharmonic on $\mathbb{C}^N.$ If $(f, D)$ is not constant, then for any holomorphic branch $(g, V)$ of $(f, D)$ in $\mathrm{Reg}(\Omega),$ it holds that $\psi(g) < 0$ on $V.$ Indeed, suppose there is a holomorphic branch $(h, U)$ of $(f, D)$ in $\mathrm{Reg}(\Omega)$ such that $\psi(h(p))=0$ for some $p \in U.$ By Theorem \ref{thm33},  $\psi(g) \equiv 0$ for every holomorphic branch $(g, V)$ of $(f, D)$.  Since $\psi$ is strictly plurisubharmonic,  we know $\{ \psi =0\}$ contains only trivial complex varieties. But $\psi \circ g \equiv 0$ on $V$. We conclude $g$ must be constant for every branch $(g, V).$ This contradicts with the assumption that $(f, D)$ is not constant.

\end{remark}

\begin{remark}\label{rekk}
In general, $(f, D)$ does not admit holomorphic continuation across
a  singular point (even normal singular point) of $\Omega$, as
demonstrated by  the following example.
\end{remark}

\begin{example}
Let $\Ol{\Omega}$ be the Stein space with boundary defined by
$$\Ol{\Omega}=\{W=(w_1, w_2, w_3) \in \mathbb{C}^3:\sum_{j=1}^3 |w_j|^2 \leq 1, w_2^2=2w_1 w_3  \}.$$
Let $\pi: \Ol{\mathbb{B}^2} \rightarrow \Ol{\Omega}$ be given by
$\pi(z_1, z_2)=(z_1^2, \sqrt{2}z_1z_2, z_2^2).$ Note $\pi$ is $2$ to
1 covering from $\mathbb{B}^2 \setminus \{ 0\}$ to $\Omega \setminus
\{ 0\}$ and $\pi(0)=0.$ Fix $p_0=(\frac{1}{2}, \frac{\sqrt{2}}{2},
\frac{1}{2}) \in M := \partial \Omega.$ Let $D$ be a small open
piece of $M$ containing $p_0$ and $(f, D)$ a CR mapping element
given by $f(W)=(\sqrt{w_1}, \sqrt{w_3}).$ Here $\sqrt{w}=\sqrt{|w|}
e^{i\frac{\theta}{2}}$ for $w=|w|e^{i \theta}$ with $-\pi< \theta <
\pi.$ Notice that $f$ maps $D$ into $\partial{\mathbb B}^2$. Since
$M$ is spherical,  $(f, D)$ admits unrestricted bounded CR
continuation in $M$ by a result of Burns-Shnider \cite{BSc}.  By Theorem
\ref{thm33}, $(f, D)$ admits unrestricted bounded holomorphic
continuation in $\mathrm{Reg}(\Ol{\Omega}).$ We claim, however, it
does not always admit holomorphic continuation across the singular
point $0$.  Set $\gamma_1$ be a curve in $\Ol{\Omega}$ such that
$\gamma_1(0) =p_0$ and $\gamma_1(1)=(\epsilon, 0, 0)$ for some small
$\epsilon > 0$ such that $\gamma_1$ never intersects $\{ w_1=0 \}$
and let the curve $\gamma_2$ in $\Ol{\Omega}$ be given by
$\gamma_2(t)=(\epsilon t, 0, 0), 0 \leq t \leq 1.$  Write
$\gamma=\gamma_1 + (-\gamma_2).$ Note $\gamma([0, 1))$ does not
pass through $\{w_1=0\}$. We know for every $0 \leq t_0 < 1,$  if we write $(h, V)$ for the branch we
obtain at $\gamma(t_0)$ with $h=(h_1, h_2)$ on $V$, then $h_1$ equals either $[-\sqrt{w_1}]_{\gamma(t_0)}$
or $[\sqrt{w_1}]_{\gamma(t_0)}.$ Without loss of generality, we
assume it is the latter when $t_0 (\neq 1)$ is close to $1$. Suppose we can extend $(f, D)$
holomorphically along
 $\gamma$ to get a holomorphic branch $(g, U)$ at $\gamma(1)=0$ (In particular $0 \in U$). Write $g=(g_1, g_2)$ on $U$. Then $g_1=\sqrt{w_1}$
near $(\epsilon_0, 0, 0)$ for a sufficiently small $\epsilon_0$. But this is impossible as we can find a loop $\sigma$ in $U$ given by $\sigma(t)=(\epsilon_0 e^{2 \pi it}, 0, 0)\in \overline{\Omega}, 0 \leq t \leq 1,$ so that we get a different branch when applying  holomorphic continuation to $\sqrt{w_1}$ along $\sigma.$
\end{example}

The proof of Theorem \ref{thm33} is  long and we will split it into
several steps, to be established in the following two subsections
$\S 3.1$ and $\S 3.2.$

Before proceeding to the proof, we first fix a Morse
plurisubharmonic defining function $\rho$ of $\overline{\Omega}.$
More precisely, we choose a bounded plurisubharmonic exhaustion
function $\rho: \overline{\Omega} \rightarrow [-\infty, 0]$ of
$\Omega$ such that $\rho \equiv 0$ on $M, \rho < 0 $ in $\Omega$ and
$\rho(z) =-\infty$ if and only if $z$ is a singular point of
$\Omega.$ In addition, $d \rho|_{M} \neq 0$ and $\rho$ is smooth
strongly plurisubharmonic on $\mathrm{Reg}(\Omega)$. Moreover,
$\rho$ has only finitely many critical points in
$\mathrm{Reg}(\Omega)$ and they are all non-degenerate. The
existence of such a $\rho$ is guaranted by the assumption on
$\overline{\Omega}$ and Morse function theory (The local existence
of such a function near a singular point can be found in Milnor
\cite{Mi}. Away from singular points, we refer to the book of
Forstneri\v{c} \cite{Fo} for such a construction. Then one applies the
Morse approximation to get our $\rho$.).

We then fix a Riemannian metric $d s^2$ over $\mathrm{Reg}(\overline{\Omega})$ which induces a distance function $\hat{d}(x, y)$ for
$x, y \in \mathrm{Reg}(\overline{\Omega}).$
Write $X_{\rho}$ for the dual vector field of $-d \rho$ with respect to $d s^2$ over the places in $\mathrm{Reg}(\overline{\Omega})$ away from the critical points of $\rho.$


\subsection{Proof of Theorem \ref{thm33}: Part I}

This step is the same for the two cases (1), (2) in Theorem \ref{thm33} and it aims to prove $(f, D)$ admits unrestricted holomorphic continuation in a tube neighborhood of $M$. We emphasize that the boundedness assumption in case (2) is not used in $\S 3.1$ to derive the extension.  We choose three finite open convex cover $\{W_j^{(k)} \}_{j=1}^m, k=1, 2, 3$ of $(M, d s^2|_M)$ with $W_j^{(1)} \ssubset W_j^{(2)} \ssubset W_j^{(3)}$ for each $j$. Moreover, we make $W_j^{(3)}$ sufficiently small for each $j$ so that a neighborhood of $\overline{W_j^{(3)}}$ on $M$ is CR diffeomorphic to a strongly pseudoconvex hypersurface in $\mathbb{C}^n.$ Write $D_j$ for the union of all smooth holomorphic disks attached to $W_j^{(3)}$ which can be deformed through a continuation family of disks to points in $W_j^{(3)}.$ For $0 < \epsilon_1 \ll 1$ and  $1 \leq k \leq 3,$ we let $\widehat{W}_{j, \epsilon_1}^{(k)}$ be the open subset of $\Ol{\Omega}$ obtained by flowing each point $p \in W_j^{(k)}$ along $X_{\rho}$ with time $0\leq  t < \epsilon_1.$ Note we can find an $\epsilon_1 >0 $ sufficiently small such that $\widehat{W}_{j, \epsilon_1}^{(2)} \ssubset (D_j \cup W_j^{(3)})$ for each $j,$ and that $\widehat{W}_{j, \epsilon_1}^{(2)}$ is topological trivial (recall $W_j^{(2)}$ is chosen to be convex). Fix such an $\epsilon_1.$ Write for $r_2 < r_1 \leq 0,$
$$\overline{\Omega}_{r_2, r_1}=\{ p \in \overline{\Omega}: r_2 < \rho \leq r_1 \}.$$
We emphasize that $\overline{\Omega}_{r_2, r_1}$ only contains  its
outer boundary but not its inner boundary.

 Let $0 < \epsilon_2 \ll \epsilon_1$ be small
enough such that $\overline{\Omega}_{-\epsilon_2, 0} \subset
\cup_{j=1}^m \widehat{W}_{j, \epsilon_1}^{(2)}.$ Define
$J_{\epsilon_2}: \Ol{\Omega}_{-\epsilon_2, 0} \rightarrow M$ for the
retract of $\Ol{\Omega}_{-\epsilon_2, 0}$ to $M$ which maps
 every point $p$ in $\overline{\Omega}_{-\epsilon_2, 0}$ through the orbit of $X_{\rho}$ to the corresponding
  point on $M$. Note $J_{\epsilon_2}$ is a smooth map for small $\epsilon_2$. By the Lewy-Baouendi-Treves
  theorem,
   we see that every continuous CR function $h$ on $W_j^{(3)}$ extends to a holomorphic function in $\widehat{W}_{j,
    \epsilon_{1}}^{(2)}$ that is continuous up to $W_j^{(2)}$.  Let $\gamma: [0, 1] \rightarrow
    \Ol{\Omega}_{-\epsilon_2, 0}$ be a curve.
    There is a corresponding curve
     $\hat{\gamma}:=J_{\epsilon_2} \circ \gamma$ on $M$. By making $\epsilon_1, \epsilon_2$ sufficiently small, we note from the definition of $J_{\epsilon_2},$  that
     $\gamma(t) \subset \widehat{W}_{j, \epsilon_1}^{(k)}$ for some $1 \leq j \leq
     m,\
      1 \leq k \leq 3$
      if and only if $\hat{\gamma}(t) \subset W_{j}^{(k)}.$
We next prove the following lemma.

\begin{lemma}\label{lemmaepx}Let $(f, D)$ and $\psi$ be as in the assumption of Theorem \ref{thm33}.
 Fix a curve $\gamma: [0, 1] \rightarrow \Ol{\Omega}_{-\epsilon_2, 0}$ with $\gamma(0) =p_0 \in D$ and
 let $\hat{\gamma}$
be as above.
\begin{enumerate}[(1)]
\item We can find $\{ W_{j_l}^{(2)}\}_{l=0}^{k}$ and $\{ \widehat{W}_{j_l, \epsilon_1}^{(2)} \}_{l=0}^{k}$
 with $0 \leq j_l \leq m$, together with CR function elements $\{ (f_l, W_{j_l}^{(2)}) \}_{l=0}^k$ on $M$
 and holomorphic function elements $\{(g_l, \widehat{W}_{j_l, \epsilon_1}^{(2)}) \}_{l=0}^k$ in
  $\mathrm{Reg}(\overline{\Omega})$ such that the following hold:
\begin{enumerate}[(a)]
\item The point $p_0 \in W_{j_0}^{(2)}$ and $f_0=f$ in a neighborhood of $p_0$ on $M$.

\item There is  a partition $0=\delta_0 < \delta_1 < \cdots < \delta_{k+1}=1$ of $[0, 1]$ such that $\gamma([\delta_l , \delta_{l+1}]) \subset \widehat{W}_{j_l, \epsilon_1}^{(2)}, \hat{\gamma}([\delta_l, \delta_{l+1}]) \subset W_{j_l}^{(2)}$ for $0 \leq l \leq k.$

\item $f_l=g_l$ on $W_{j_l}^{(2)}$ for $0 \leq l \leq k.$

\item $f_l =f_{l+1}$ on $W_{j_l}^{(2)}\cap W_{j_{l+1}}^{(2)}$ and $g_l=g_{l+1}$ on $\widehat{W}_{j_l, \epsilon_1}^{(2)} \cap \widehat{W}_{j_{l+1}, \epsilon_1}^{(2)}$ for $0 \leq l \leq k-1.$
\end{enumerate}

Consequently, $\{ (g_l, \widehat{W}_{j_l, \epsilon_1}^{(2)})\}_{l=0}^k$ (resp. $\{(f_l,  W_{j_l}^{(2)}) \}_{l=0}^k$) induces a holomorphic (resp. CR) continuation of $(f, D)$ along $\gamma$ (resp. $\hat{\gamma}$).

\item It holds that $\psi \circ g_l \leq 0$ on $\widehat{W}_{j_l, \epsilon_1}^{(2)}$ for all $0 \leq l \leq k$.  Moreover, assume
there is some $0 \leq l_0 \leq k$ and a point $q \in \widehat{W}_{j_{l_0}, \epsilon_1}^{(2)} \setminus M$ such that $\psi \circ g_{l_0}(q) =0.$ Then $\psi(g) \equiv 0$ for any holomorphic branch $(g, V)$ of $(f, D)$ in $\Ol{\Omega}_{-\epsilon_2, 0}.$


\end{enumerate}
\end{lemma}

{\bf Proof of Lemma \ref{lemmaepx}:} By the uniform continuity of
$\hat{\gamma}$ on $[0, 1]$ and Lebesgue's number lemma, we can find some
$\epsilon
>0$ such that for any sub-interval $I^*$ of $[0,1]$ with length
bounded by $\epsilon$, there exists some $1 \leq j(I^*) \leq m$ satisfying
$\hat{\gamma}(I^*) \subset W_{j(I^*)}^{(1)}$.
Note $p_0 =\gamma(0)=\hat{\gamma}(0)$ is contained in
$W_{j_0}^{(1)}$ for some $1 \leq j_0 \leq m.$ Set $\delta_0=0$ and
let $\delta_1 \in (0, 1]$ be the (unique) number (if exists) such
that $\hat{\gamma}([0, \delta_1)) \subset W_{j_0}^{(1)}$ and
$\hat{\gamma}(\delta_1) \not\in W_{j_0}^{(1)}$. Note here we
choose $j_0$ such that $\delta_1$ takes the largest value and thus we
must have $\delta_1 \geq \epsilon$ if it exists. If such a number
$\delta_1$ does not exist, this means $\hat{\gamma}([0,1]) \subset
W_{j_0}^{(1)}$ and consequently $\gamma ([0, 1]) \subset
\widehat{W}_{j_0, \epsilon_1}^{(1)} \cap \Ol{\Omega}_{-\epsilon_2,
0}.$ Note first by the unrestricted CR continuation assumption, the
germ of $f$ at $p_0$ extends to a CR function $f_0$ on
$W_{j_0}^{(3)}$ as $W_{j_0}^{(3)}$ is simply connected. Secondly
$(f_0, W_{j_0}^{(2)})$ can be extended to a holomorphic function
element $(g_0, \widehat{W}_{j_0, \epsilon_1}^{(2)})$. Then $(g_0,
\widehat{W}_{j_0, \epsilon_1}^{(2)})$ induces a holomorphic
continuation of $(f, D)$ along $\gamma$ and the first part of lemma
is established.

Now assume such a $\delta_1$ exists. First, as above the germ of $f$ at $p_0$ extends to a CR function $f_0$ on $W_{j_0}^{(3)}$ and $f_0$ extends to a holomorphic function element $(g_0, \widehat{W}_{j_0, \epsilon_1}^{(2)})$(Thus $(c)$ holds for $l=0$).  We then look at $\hat{\gamma}(\delta_1)$. Note $\hat{\gamma}(\delta_1) \in W_{j_0}^{(2)} \subset W_{j_0}^{(3)}$ and there exists some $1 \leq j_1 \leq m$ such that
$\hat{\gamma}(\delta_1) \in W_{j_1}^{(1)}.$ By the same reason as above, the germ of $f_0$ at $\hat{\gamma}(\delta_1)$  extends to a CR function $f_1$ on $W_{j_1}^{(3)}$ and $f_1$ extends to a holomorphic function element $(g_1, \widehat{W}_{j_1, \epsilon_1}^{(2)})$(Thus (c) holds for $l=1$). Note $f_1$ and $f_0$ coincide near $\hat{\gamma}(\delta_1)$. Moreover, $W_{j_0}^{(3)} \cap W_{j_1}^{(3)}$ and $\widehat{W}_{j_0, \epsilon_1}^{(2)} \cap \widehat{W}_{j_1, \epsilon_1}^{(2)}$ are simply connected by our construction. We conclude that $f_0=f_1$ on $W_{j_0}^{(3)} \cap W_{j_1}^{(3)}$ and $g_1=g_0$ in $\widehat{W}_{j_0, \epsilon_1}^{(2)} \cap \widehat{W}_{j_1, \epsilon_1}^{(2)}$(Thus $(d)$ holds for $l=0$).

Then we pick the (unique) number $\delta_2$ (if exists) such that
$\hat{\gamma}([\delta_1, \delta_2)) \subset W_{j_1}^{(1)}$ and
$\hat{\gamma}(\delta_2) \not\in W_{j_1}^{(1)}$. Note again we
choose $j_1$ such that $\delta_2$ takes the largest possible value
if exits and thus we must have $\delta_2-\delta_1 \geq \epsilon$. And we
run the same procedure as above to obtain $1 \leq j_2 \leq m$ such
that $\hat{\gamma}(\delta_2) \in W_{j_2}^{(1)}$, along with a CR
function element  $(f_2, W_{j_2}^{(3)})$ and a holomorphic function
element $(g_2, \widehat{W}_{j_2, \epsilon_1}^{(2)})$ such that
$f_1=f_2$ on $W_{j_1}^{(3)} \cap W_{j_2}^{(3)}$, and $g_1=g_2$ in
$\widehat{W}_{j_1, \epsilon_1}^{(2)} \cap \widehat{W}_{j_2,
\epsilon_1}^{(2)}.$

By repeating the above procedure for at most $[\frac{1}{\epsilon}]+1$ times, we arrive at
some positive number $\delta_{k}$ such that $\delta_{k} + \epsilon \geq 1.$  More precisely,
 we obtain a partition $0 =\delta_0 < \delta_1 < \cdots < \delta_{k} < \delta_{k+1}=1$ of $[0, 1]$
 and a collection of integers $1 \leq j_1, \cdots, j_k \leq m$ such that
  $\hat{\gamma}([\delta_l, \delta_{l+1})) \in W_{j_l}^{(1)}$ and $\hat{\gamma}(\delta_{l+1})
  \not \in W_{j_l}^{(1)}$ for $1 \leq l \leq k$
   (In particular, $\hat{\gamma}([\delta_l, \delta_{l+1}]) \in W_{j_l}^{(2)}$ for all $k$).
    Moreover, there are a collection of CR function elements $\{ (f_l, W_{j_l}^{(3)})\}_{l=0}^k$ and
    holomorphic function elements $\{ (g_l, \widehat{W}_{j_l, \epsilon_1}^{(2)})\}_{l=1}^k$ satisfying
     condition $(c)$ and $(d).$ This proves Part (1) of the lemma.

To prove Part (2), we note every $p \in \widehat{W}_{j_l, \epsilon_1}^{(2)} \setminus M,$ there exists a smooth holomorphic disk $\Delta_p$ attached to $W_{j_l}^{(3)}$ such that $p \in \Delta_p$. Moreover, $\psi \circ g_l$ is subharmonic in $\Delta_p$ and is continuous up to $\partial \Delta_p,$ and it agrees with $\psi \circ f_l$ on $\partial \Delta_p$. By assumption $(f, D)$ admits unrestricted CR continuation with $\psi-$estimate in $M$ and thus in particular $\psi \circ f_l \leq 0$ on $ \partial \Delta_p  \subset W_{j_l}^{(3)}$, we conclude by the maximum principle that $\psi \circ g_l (p) \leq 0.$ As $p$ is arbitrary, we have $\psi \circ g_l \leq 0$ on $\widehat{W}_{j_l, \epsilon_1}^{(2)}$ for all $0 \leq l \leq k.$ 
Now suppose there is some $0 \leq l_0 \leq k$ and a point $q \in \widehat{W}_{j_{l_0}, \epsilon_1}^{(2)} \setminus M$ such that $\psi \circ g_{l_0}(q) =0.$ This means $\psi \circ g_{l_0}$ acheives its maximum at an interior point $q$.  Since
$\psi \circ g_{l_0}$ is subharmonic in $\widehat{W}_{j_{l_0}, \epsilon_1}^{(2)},$ we conclude that $\psi \circ g_{l_0} \equiv 0$ in $\widehat{W}_{j_{l_0}, \epsilon_1}^{(2)}.$ As $g_l=g_{l+1}$ on $\widehat{W}_{j_l, \epsilon_1}^{(2)} \cap \widehat{W}_{j_{l+1}, \epsilon_1}^{(2)}, $ we have each $\psi \circ g_l$ attains its maximum at an interior point and thus is constant. In particular, $\psi \circ g_0 \equiv 0$ in $\widehat{W}_{j_0, \epsilon_1}^{(2)}$ and $\psi \circ f_0=0$ in $W_{j_0}^{(2)}.$ Now let $(g, U)$ be a  holomorphic branch of $(f, D)$ in $\Ol{\Omega}_{-\epsilon_2, 0}.$
Then there is a path $\sigma$ in $\Ol{\Omega}_{-\epsilon_2, 0}$ connecting $p_0$ to some point $p_1 \in U$ such that $(g, U)$ is obtained by holomorphic continuation of $(f, D)$ along $\sigma.$
By part (1) of Lemma \ref{lemmaepx},  writing $\hat{\sigma}:=J_{\epsilon_1} \circ \sigma,$ there exist holomorphic function elements $\{(\widetilde{g}_l, \widehat{W}_{i_l, \epsilon_1}^{(2)}) \}_{l=0}^{\nu}$ and CR function elements $\{(\widetilde{f}_l, W_{i_l, \epsilon_1}^{(2)}) \}_{l=0}^{\nu}$ that satisfy the conditions (a), (b), (c), and (d), and induce a holomorphic continuation and a CR continution of $(f, D)$ along $\sigma$ and $\hat{\sigma}.$  By uniqueness of holomorphic continuation, we have $g=\hat{g}_{\nu}$ near $q$ in $U \cap \widehat{W}_{i_{\nu}, \epsilon_1}^{(2)}$.  Note we have $\widetilde{f}_0=f_0$ near $p_0$ and $\widetilde{g}_0=g_0$ near $p_0.$
Then we have $\psi \circ \widetilde{g}_0 \equiv 0$ on $\widehat{W}_{i_0, \epsilon_1}^{(2)}$. Applying the maximum principle for subharmonic functions as above, we obtain that $\psi \circ \widetilde{g}_l \equiv 0$ in $\widehat{W}_{i_l, \epsilon_1}^{(2)}$ for every $0 \leq l \leq \nu.$ In particular $\psi \circ g \equiv 0$ on $U$. This proves Lemma
\ref{lemmaepx}.  \endpf

\bigskip


By Lemma \ref{lemmaepx}, we immediately have the following corollary.

\begin{corollary}\label{coromegae} Let $(f,D)$ be as in Theorem
\ref{thm33}.  The function element $(f, D)$ admits unrestricted
holomorphic continuation in $\Ol{\Omega}_{-\epsilon_2, 0}$ with
$\psi-$estimate (In the dimension two case (2), $(f, D)$ admits
unrestricted bounded holomorphic continuation in
$\Ol{\Omega}_{-\epsilon_2, 0}$ with $\psi-$estimate).  Moreover, if
there is a holomorphic branch $(h, U)$ of $(f, D)$ in
$\Ol{\Omega}_{-\epsilon_2, 0}$ such that $\psi (h(p))=0$ at some
point $p$ in $U \setminus M$. Then $\psi(g) \equiv 0$ for every
holomorphic branch $(g, V)$ of $(f, D)$ in
$\Ol{\Omega}_{-\epsilon_2, 0}$.
\end{corollary}

\subsection{Proof of Theorem \ref{thm33}: Part II}

In this step,  we finish the proof of Therem \ref{thm33}. We
emphasize that the treatment for the two cases are  different. We
will apply the method of continuous family of holomorphic curves
which is a  typical machinery in the study of holomorphic
continuation problem.
 The use of the Morse function theory to study the holomorphic
 continuation of multiple-valued holomorphic functions near the
 boundary to the interior of the pseudo-convex domain
  first appeared in $\S 5$ of the paper
 by Huang-Ji \cite{HJ}. In the  paper by Merkel-Porten \cite{MP}, they   employed the Morse
 function theory to re-investigate  the Hartogs extension theorem.
 In our argument here, besides the Morse function theory, the
Phragm\'en-Lindel\"of principle will also play a fundamental role.

To start with, we make the following definition. Let $(f, D)$ be as
in Theorem \ref{thm33}, i.e., it satisfies the assumption in (1) in
the three or higher dimensional case, and satisfies the assumption
in (2) in the dimension two case. Let $\widetilde{\Omega}$ be a
connected open subset of $\mathrm{Reg}(\Ol{\Omega})$ with $D \subset \widetilde{\Omega}.$ We say
$\widetilde{\Omega}$ has the extendable property if, in the three or
higher dimensional case, $(f, D)$ admits unrestricted holomorphic
continuation with $\psi-$estimate in $\widetilde{\Omega}$; and in
the dimension two case, $(f, D)$ admits unrestricted bounded
holomorphic continuation in $\widetilde{\Omega}$ with
$\psi-$estimte, where $D\cap \Omega\subset \widetilde{\Omega}$. To
establish Theorem \ref{thm33}, we will need to prove
$\mathrm{Reg}(\Ol{\Omega})$ has the extendable property. Now set
$$A=\{a < 0: \text{The set}~\Ol{\Omega}_{a, 0}~\text{has the extendable property}\}.$$
We first note that $A$ is not empty as $-\epsilon_2 \in A.$ Set $b =\mathrm{inf}(A) < 0.$ Note if $b=-\infty,$ then Theorem \ref{thm33} is established.  Indeed note every curve $\gamma: [0, 1] \rightarrow \mathrm{Reg}(\Ol{\Omega})$ must be contained in $\Ol{\Omega}_{c, 0}$ for some $c < 0$ and thus $\mathrm{Reg}(\Ol{\Omega})$ has the extendable property.   We will therefore assume $b > -\infty$ in the following context. Note it follows from the definition of $b$ that $\Ol{\Omega}_{b, 0}$ has the extendable property. Writing $\mathrm{inf} (\rho):=\mathrm{inf}\{\rho(z): z \in  \Ol{\Omega} \setminus \mathrm{sing}(\Omega) \},$ we have two cases:

\medskip
{\bf Case I:} We first consider the case if $b =\mathrm{\inf}(\rho).$ In this case, since $b > -\infty,$ we must have that $\mathrm{sing}(\Omega)= \emptyset.$ Write $\rho^{-1}(c)$ for the level set $\{z \in \Ol{\Omega}: \rho(z)=c \}.$ Note in this case we have $\Ol{\Omega}_{b, 0}=\mathrm{Reg}(\Ol{\Omega}) \setminus \rho^{-1}(b)$. Fix any $\hat{q} \in \rho^{-1}(b).$ By assumption $\hat{q}$ is an isolated critical point of $\rho$. Moreover, the real Hessian of $\rho$ is nondegenerate and thus strictly positive definite at $\hat{q}$. For every such $\hat{q},$ choose a small neighborhood $U_{\hat{q}}$ of $\hat{q}$ and real  coordinates $x$ with $x(\hat{q})=0$ such that $\hat{q}$ is the only critical point in $\Ol{U_{\hat{q}}}$ and  $\rho(x)=\sum_{j=1}^{2n} x_j^2 +b$ in $U_{\hat{q}}$. Then there exists a negative number $c(\hat{q}) > b$ such that $\rho^{-1}(c) \cap U_{\hat{q}} \ssubset U_{\hat{q}}$ is a sphere whenever $b < c \leq c(\hat{q})$.  Moreover, this sphere is of real dimension at least three, and thus is simply connected.  Take $c_0 >b$ to be smaller than $c(\hat{q})$ for all $\hat{q} \in \rho^{-1}(b).$ And write $V(\hat{q}, c_0):=U_{\hat{q}} \cap \{ \rho (x) < c_0 \} \ssubset U_{\hat{q}}$ for $\hat{q} \in \rho^{-1}(b).$ Let $\gamma: [0, 1] \rightarrow \mathrm{Reg}(\Ol{\Omega})$ be any curve with $\gamma(0) \in D$. We will show $(f, D)$ admits holomorphic continuation along $\gamma$ with $\psi-$estimate (and with the boundedness in dimension two case). Note this  follows immediately from the fact that $\mathrm{Reg}(\Ol{\Omega}) \setminus \rho^{-1}(b)$ has the extendable property if $\gamma$ does not intersect $\rho^{-1}(b).$  Now assume $\gamma$ intersects $\rho^{-1}(b).$ By the uniform continuity of $\rho \circ \gamma$ on $[0, 1],$ there exists some small $\delta >0$ such that whenever $\rho (\gamma(t^*))=\hat{q} \in \rho^{-1}(b)$ and $|t-t^*| \leq \delta,$ we have $\gamma(t) \in V(\hat{q}, c_0) \ssubset  U_{\hat{q}}$. Let $t_0$ be the minimum number such that $\gamma(t_0) \in \rho^{-1}(b)$ and write $q_0=\gamma(t_0).$ Let $t_2=t_0+\delta$(Or, let it be $1$ if the number does not lie in $[0, 1]$). Let $t_1$ be the largest number $t$ such that $t < t_0$ and $\gamma(t) \not \in V(q_0, c_0).$
Then $\rho(\gamma(t_1)) =c_0$ and for every $t \in (t_1, t_2], \rho(\gamma(t)) \in V(q_0, c_0) \ssubset U_{q_0}$.
By assumption, $(f, D)$ admits holomorphic continuation along $\gamma([0, t_1])$ and we get a holomorphic branch at the point $q=\gamma(t_1) \in \rho^{-1}(c_0)$. Write this branch as $(g, B).$ Let $W_0 \ssubset W_1 \ssubset U_{q_0}$ be two small simply connected tube neighborhoods of $\rho^{-1}(c_0) \cap U_{q_0}$. As $W_1  \subset \mathrm{Reg}(\Ol{\Omega}) \setminus \rho^{-1}(b)$, $(g, B)$ can be holomorphically extended along any path in $W_1$, thus it induces a well-defined holomorphic function $G$ in $W_1$.  And $(G, W_1)$ satisfies the $\psi-$estimate by assumption. Then by the Hartogs extension, $G$ extends to a holomorphic function, which we still denote by $G$, to an open subset $W^*$  that contains $q_0$. Here $W^*$ equals the union of $V(q_0, c_0)$ and $W_1$.
Note $G$ also satisfies the $\psi-$estimate on $W$ by the maximum principle.  In this way, we are able to extend $(f, D)$ holomorphically along  $\gamma$ over the interval $[0, t_2]$ with the desired estimate. In the dimension two case, the boundedness of the continuation also follows from the maximum principle.

Write $W=V(q_0, c_0) \cup W_0.$ Note $\gamma([t_1, t_2]) \subset W \ssubset W^*$. If $\gamma([t_2, 1]) \subset W$, then the proof is finished. Now assume there is some $t > t_2$ in $[0, 1]$ such that $\gamma(t) \not \in W$, we let $t_3$ be the smallest such number. We can perturb $\gamma([t_1, t_3])$ in $\Ol{W}$, with the endpoints fixed, to a new curve $\widetilde{\gamma}([t_1, t_3])$ such that it avoids $q_0$ and gives the same branch of function at the point $\gamma(t_3).$ In this way, $\gamma([0, t_1])+ \widetilde{\gamma}([t_1, t_3])$  lies in $\mathrm{Reg}(\Ol{\Omega}) \setminus \rho^{-1}(b).$
If $\gamma$ intersects $\rho^{-1}(b)$ again after $t_3,$ we repeat the above Hartogs extension argument. After applying this argument for at most $[\frac{1}{\delta}]+1$ times (recall $|t_2-t_0| =\delta$), we will finally arrive at $\gamma(1).$ Thus $(f, D)$ admits holomorphic continuation along $\gamma$ with the desired estimate. As $\gamma$ is arbitrary, this implies $\mathrm{Reg}(\Ol{\Omega})$ has the extendable property. This finishes the proof in case I.

\bigskip

Now it remains to consider the case $b > \mathrm{inf}(\rho).$ We will prove this indeed cannot happen.

\medskip
{\bf Case II:} We first assume $b > \mathrm{inf}(\rho)$ and $M_{b}:=\rho^{-1}(b)$ has no critical points
 of $\rho.$ In this case, we note $M_{b'}$ is a smooth strongly pseudoconvex hypersurface in $\Omega$ for $b'\approx b$.
  Recall by the definition of $b$, $\Ol{\Omega}_{b, 0}$ has the extendable property.
 Then we apply the same
argument in $\S3.1$  with $b'$ being sufficiently close but greater
than $b$ to obtain a small $\epsilon
>0$ such that $(f, D)$ admits unrestricted holomorphic extension in
$\Ol{\Omega}_{b-\epsilon, 0}$ with $\psi-$estimate. We thus conclude
$\Ol{\Omega}_{b-\epsilon, 0}$ has the extendable property. This is,
however, a contradiction to the definition of $b$.

\medskip
{\bf Case III:} If $b > \mathrm{inf}(\rho)$ and $\rho$ has critical
points on $M_b.$ Let $p \in \rho^{-1}(b)$ be a critical point of
$\rho$. Then choose a neighborhood $U_{p}$ of $p$ such that $p$ is
the only critical point of $\rho$ in $U_p.$ And choose certain
holomorphic coordinates $z$ on $U_p$ such that $z(p)=0$ and $\rho$ takes the
following form near $p$:
\begin{equation}
\rho=||z||^2 + 2 \mathrm{Re}\sum_{j=1}^n \lambda_j z_j^2 + O(||z||^3)+b.
\end{equation}
 Here we have $0 \leq  \lambda_j < \infty$ and $\lambda_j \neq \frac{1}{2}$ by the non-degeneracy assumption. Recall that the $z_j-$direction is called elliptic if $0 \leq \lambda_j < \frac{1}{2},$ and hyperbolic if $\lambda_j > \frac{1}{2}.$ Also, in some smooth coordinates $x$ on $U_{p}$ with $x(p)=0,$ we have
$$\rho(x)=\sum_{j=1}^m x_j^2-\sum_{m+1}^{2n} x_j^2+b. $$
By the plurisubharmonicity, we have $m \geq n.$ For any small number $\epsilon >0,$ set $V_{\epsilon}:=\{ q \in U_{p}: ||x(q)|| < \epsilon\} \ssubset U_{p}.$ One can directly verify that $V_{\epsilon} \cap \Ol{\Omega}_{b, 0}$ is connected for small $\epsilon >0.$ We will need the following crucial lemma:
\begin{lemma}\label{lemmaqvepsi}
For sufficiently small $\epsilon_3 >0,$  it holds for every $q \in V_{\epsilon_3}
 \cap \Ol{\Omega}_{b, 0},$  that if $[g]_q$ is a holomorphic branch of $(f , D)$ obtained by holomorphic
  continuation along a curve in $\Ol{\Omega}_{b, 0}$, then $[g]_q$ extends to a single-valued holomorphic function in $V_{\epsilon_3}$ with $\psi-$estimate (and with the boundedness in the dimension two case).
\end{lemma}

\begin{proof}
We choose holomorphic coordinates $z$ in a small neighborhood
$U=U_p$ of $p$ mentioned above (in particular (7) holds and $p$ is
the only critical point of $\rho$ in $U$). It holds that for a
sufficiently small $\delta$ and a small $||z||$ that
\begin{equation}\label{eqnrhoest}
 \rho \geq (1-\delta)||z||^2 + 2 \mathrm{Re}\sum_{j=1}^n \lambda_j z_j^2 + b.
\end{equation}

{\bf (a):} We first assume $\rho$ has an elliptic direction at $z(p)=0$, say the $z_{1}-$direction(i.e., $\lambda_1 < \frac{1}{2}$). Write $\Delta$ for the unit disk in $\mathbb{C}.$ Fix  small numbers $0 <\widetilde{\epsilon} \ll 1$ and  $0 < \eta \ll 1.$ We define for $0 < t < \widetilde{\epsilon}, $ a continuous family (parametrized by $t$) of holomorphic disks with boundary $\phi_t(\xi): \Ol{\Delta} \rightarrow U_0 \ssubset U_p$ given by $(\eta \xi, \eta t, 0, \cdots, 0)$.
By (\ref{eqnrhoest}) we have for $\xi \in \Ol{\Delta},$
\begin{equation}\label{eqnrhophi}
\rho \circ \phi_t (\xi) \geq (1-\delta)\eta^2|\xi|^2 + \lambda_1 \eta^2( \xi^2 + \Ol{\xi}^2) +(1- \delta + 2 \lambda_2) \eta^2 t^2  +b.
\end{equation}
Choosing $\delta$ small enough such that $\frac{\lambda_1}{1-\delta} < \frac{1}{2},$ we have
$$\rho \circ \phi_t (\xi) \geq  b + (1-\delta + 2 \lambda_2) \eta^2 t^2 > b. $$
Thus we have $\phi_t (\Ol{\Delta}) \ssubset \Ol{\Omega}_{b, 0}$. Now for $\xi \in \partial \Delta,$ it follows from (\ref{eqnrhophi}) again that
\begin{equation}\label{eqnrhophitxi}
\rho (\phi_t (\xi)) \geq (1-\delta-2 \lambda_1) \eta^2 +b.
\end{equation}
In the following context, we set $d(z, w)=\mathrm{max}\{|z_j-w_j|: 1
\leq j \leq n\}$ for $z=(z_1, \cdots, z_n)$ and $w=(w_1, \cdots,
w_n)$ in $\mathbb{C}^n.$ Set $d(X, Y)=\mathrm{inf}\{d(z, w): z \in
X, w \in Y\}$ for two subsets $X, Y$ of $\mathbb{C}^n.$ Write
$\mathbb{P}(z_0, r) \subset \mathbb{C}^n, r > 0, $ for the polydisk
$\{z \in \mathbb{C}^n: d(z, z_0) < r \}.$

Note (\ref{eqnrhophitxi}) implies there exists a positive number $A_1$ independent of $t$ such that
\begin{equation}\label{eqnradiusa}
d(\phi_t(\partial \Delta), \partial (\Ol{\Omega}_{b, 0} \cap U)) \geq A_1~\text{for all}~t.
\end{equation}
Note there also exists a positive number $A_2$ independent of $t$  such that
\begin{equation}\label{eqnradiusda}
d(\phi_t(\Ol{\Delta}), \partial U) \geq A_2~\text{for all}~t.
\end{equation}
Set $A=\mathrm{min}\{ A_1, A_2 \}.$ Now pick $\epsilon >0$
sufficiently small such that $V_{\epsilon} \subset U$ and  $d(z,
0) < \frac{A}{2}$ whenever $z \in V_{\epsilon}.$ Now fix any branch
$[g]_q$ with $q \in V_{\epsilon} \cap \Ol{\Omega}_{b, 0}$ as in the
assumption of Lemma \ref{lemmaqvepsi}. As $V_{\epsilon} \cap
\Ol{\Omega}_{b, 0}$ is connected, we can first extend $[g]_q$
holomorphically along certain curve in $V_{\epsilon} \cap
\Ol{\Omega}_{b, 0}$ to obtain a new branch $[h]_{\hat{q}}$ with
$\hat{q}= \phi_t(0)=(0, \eta t, 0, \cdots, 0) \in V_{\epsilon} \cap
\phi_t (\Ol{\Delta})$ for sufficiently small $0 <t
<\widetilde{\epsilon}$.

We also note that
$$\mathrm{d}({\bf 0}, \phi_t(0))=\eta t \rightarrow 0~\text{as}~t \rightarrow 0.$$
As for each fixed small $t,$ $\phi_t (\Ol{\Delta}) \ssubset \Ol{\Omega}_{b, 0},$ we can extend $[h]_{\hat{q}}$ holomorphically along any curve inside a small neighborhood of $\phi_t (\Ol{\Delta}).$ Since $\phi_t (\Ol{\Delta})$ is simply connected, $[h]_{\hat{q}}$ extends to a well-defined holomorphic function, which we still denote by $h$, in this small
neighborhood of $\phi_t (\Ol{\Delta})$.
Moreover, if $q_0 \in \phi_t(\partial \Delta)$, then by (\ref{eqnradiusa}) and (\ref{eqnradiusda}), we can extends $[h]_{q_0}$ to a holomorphic function in the polydisk $\mathbb{P}(q_0, A)$. By the continuity principle, we conclude that for any $w \in \phi_t (\Ol{\Delta}),$ the Taylor expansion (in $z-$coordinates) of  $[h]_w$ about $p$ converges in $\mathbb{P}(w, A).$ In particular, $[h]_{\hat{q}}$ (Recall $\hat{q}=\phi_t(0)$) extends to a holomorphic function, still called $h,$ in $\mathbb{P}(\hat{q}, A)$ and note for sufficiently small $t, \mathbb{P}(\hat{q}, A)$ contains $V_{\epsilon}.$ By the uniqueness of holomorphic maps, we have $[g]_q=[h]_q.$ In this way, we extends $[g]_q$ to a holomorphic map $h$ in $V_{\epsilon}.$

To obtain the $\psi-$estimate on $h$, we need to shrink $V_{\epsilon}.$ For  small $\lambda, \mu >0,$ define a complex $(n-1)-$parameter family(parametrized by $\tau \in \mathbb{C}^{n-1}$) of holomorphic disks  $\varphi_{\tau}(\xi): \Delta \rightarrow U$ given by $\varphi_{\tau}(\xi)=(\lambda \xi, \tau)$ with $|| \tau || < \mu.$ Note there exists $c$ that only depends on $\epsilon$ such that if we choose $\lambda, \mu < c$, then $\varphi_{\tau}(\Ol{\Delta}) \subset V_{\epsilon}$ for all $||\tau|| < \mu.$ Now we fix $0 < \mu < \lambda < c$ such that $\varphi_{\tau}(\partial \Delta) \subset \Ol{\Omega}_{b, 0}$
for all $||\tau|| < \mu$(the existence of such $\mu, \lambda $ is due to (\ref{eqnrhoest}) and the fact that $\lambda_1 < \frac{1}{2}$). Then we note there is a small $0 < \epsilon_3 < \epsilon$ which only depends on $\lambda$ and $\mu$ such that $V_{\epsilon_3} \subset \cup_{||\tau || < \mu} \varphi_{\tau}(\Ol{\Delta}).$ We claim $V_{\epsilon_3}$ is the desired region in Lemma \ref{lemmaqvepsi}. Indeed, for any $q \in V_{\epsilon_3} \cap \Ol{\Omega}_{b, 0},$ if $[g]_q$ is a holomorphic branch of $(f, D)$ in $\Ol{\Omega}_{b, 0},$ then by the above argument, $[g]_q$ extends to a holomorphic function $h$ in $V_{\epsilon}.$ Furthermore, as $V_{\epsilon} \cap \Ol{\Omega}_{b, 0}$ is connected, for each $\tau$ we can find some path in $V_{\epsilon} \cap \Ol{\Omega}_{b, 0}$ connecting $q$ to any point $q_0 \in \varphi_{\tau}(\partial \Delta).$ This implies $[h]_{q_0}$ is a branch of $(f, D)$ obtained by holomorphic continuation along a curve in $\Ol{\Omega}_{b, 0}.$ By the assumption on $\Ol{\Omega}_{b, 0}$, we conclude $[h]_{q_0}$ satisfies the $\psi-$estimate for every $q_0 \in \varphi_{\tau}(\partial \Delta)$. Now by the maximum principle for subharmonic functions, $h$ also satisfies the $\psi-$estimate on $\varphi_{\tau}(\Ol{\Delta})$ for all $\tau.$ In particular, $h$ satisfies the $\psi-$estimate in $V_{\epsilon_3}.$ The boundedness in the dimension two case follows similarly.

{\bf (b): } Next we consider the case where there are no elliptic directions (i.e., all $\lambda_j > \frac{1}{2}$).  When $n \geq 3,$ we replace the holomorphic disks $\phi_t(\xi)$ by
$$\phi_{t}(\xi)=\eta(\sqrt{\lambda_2} \xi, i \sqrt{\lambda_1} \xi, t, 0, \cdots, 0), 0 < t < \widetilde{\epsilon},$$
where $\eta$ and $\widetilde{\epsilon}$ are fixed small positive numbers.
Then by (\ref{eqnrhoest}) we have for $\xi \in \Ol{\Delta},$
$$\rho \circ \phi_t(\xi) \geq (1-\delta)(\lambda_1 +\lambda_2)\eta^2 |\xi|^2 +
(1-\delta)\eta^2t^2 +2 \lambda_3 \eta^2 t^2 +b \geq b. $$ Thus we
have $\phi_t (\Ol{\Delta}) \ssubset \Ol{\Omega}_{b, 0}$ for all $t$.
Furthermore, if $\xi \in \partial \Delta,$ it yields that
\begin{equation}\label{eqnlda12}
\rho \circ \phi_t(\xi) \geq (1-\delta)(\lambda_1 + \lambda_2) \eta^2+b.
\end{equation}

Replacing (\ref{eqnrhophitxi}) by (\ref{eqnlda12}), the same argument in {\bf (a)} yields that for sufficiently small $\epsilon,$ and for every $q \in V_{\epsilon} \cap \Ol{\Omega}_{b, 0},$ if $[g]_q$ is a holomorphic branch of $(f , D)$ in $\Ol{\Omega}_{b, 0}$, then $[g]_q$ extends to a single-valued holomorphic function $h$ in $V_{\epsilon}$. It remains to establish the desired $\psi-$estimate. We will also need to shrink $V_{\epsilon}.$ For fixed small $\lambda, \mu> 0,$ we define a complex parameter family (parametrized by $\tau \in \mathbb{C}, \chi \in \mathbb{C}^{n-2}$) of Riemann surfaces with boundary $\Ol{E}_{\tau, \chi},$ with $|\tau|, ||\chi|| < \mu,$ given by
$$\Ol{E}_{\tau, \chi}:=\{z=(z_1, z_2, \chi) \in U :\lambda_1 z_1^2 +\lambda_2 z_2^2=\tau;~~~~2 \lambda_1 |z_1|^2 + 2 \lambda_2 |z_2|^2 \leq \lambda^2 \}.$$
Note there exists $\hat{\lambda} > 0$ such that if we choose $0 < \lambda, \mu < \hat{\lambda},$ then $\Ol{E}_{\tau,\chi} \subset V_{\epsilon}$ for all $|\tau|, ||\chi|| < \mu$. Furthermore, note for $z \in \partial \Ol{E}_{\tau, \chi},$ we have
$$\rho(z) \geq (1-\delta)||z||^2 + 2 \mathrm{Re} \tau + 2 \mathrm{Re} \sum_{j=3}^n \lambda_j z_j^2 +b \geq  (1-\delta)\frac{\lambda^2}{ 2 \lambda_1 + 2 \lambda_2}-2\mu-2\mu^2(\sum_{j=3}^n \lambda_j) +b.$$
Thus we can choose $0 < \mu  \ll \lambda < \hat{\lambda}$ such that $\partial \Ol{E}_{\tau, \chi} \subset \Ol{\Omega}_{b, 0}$ for all $|\tau|, ||\chi|| < \mu.$ Fix such a pair $\lambda$ and $\mu.$
Then there is a small $0 < \epsilon_3 < \epsilon$ which only depends on $\lambda$ and $\mu$ such that $V_{\epsilon_3} \subset \cup_{|\tau |,||\chi|| < \mu} \Ol{E}_{\tau, \chi}.$ We claim $V_{\epsilon_3}$ is the desired region in Lemma \ref{lemmaqvepsi}. Indeed, for any $q \in V_{\epsilon_3} \cap \Ol{\Omega}_{b, 0},$ if $[g]_q$ is a holomorphic branch of $(f, D)$ in $\Ol{\Omega}_{b, 0},$ then by the above argument, $[g]_q$ extends to a holomorphic function $h$ in $V_{\epsilon}.$ Furthermore, as $V_{\epsilon} \cap \Ol{\Omega}_{b, 0}$ is connected, for each $\tau, \chi$ we can find some path in $V_{\epsilon} \cap \Ol{\Omega}_{b, 0}$ connecting $q$ to any point $q_0 \in \partial \Ol{E}_{\tau, \chi}.$ This implies $[h]_{q_0}$ is a branch of $(f, D)$ obtained by holomorphic continuation along a curve in $\Ol{\Omega}_{b, 0}.$ By the assumption on $\Ol{\Omega}_{b, 0}$, we conclude $[h]_{q_0}$ satisfies the $\psi-$estimate. Now by the maximum principle for subharmonic functions, $h$ also satisfies the $\psi-$estimate on $\Ol{E}_{\tau, \chi}$ for all $\tau$ and $\chi.$ In particular, $h$ satisfies the $\psi-$estimate in $V_{\epsilon_3}.$

{\bf (c):} We then consider the more difficult case when $n=2$ and
both directions are non-elliptic, i.e., $\lambda_1, \lambda_2 >
\frac{1}{2}$. Recall in this dimensional two case (2), we
additionally assume the holomorphic continuation of $(f, D)$ in
$\Ol{\Omega}_{b, 0}$ is bounded. Fix $0 < \widetilde{\epsilon} <
\eta \ll 1.$ Consider a continuous family (parametrized by $t$) of
Riemann surfaces with boundary $\Ol{E}_t, 0 < t <
\widetilde{\epsilon},$ in $U$ defined by
$$\Ol{E}_t:=\{z=(z_1, z_2) \in U :\lambda_1 z_1^2 +\lambda_2 z_2^2=t;~~~~2 \lambda_1 |z_1|^2 + 2 \lambda_2 |z_2|^2 \leq \eta^2 \}.$$
We apply a holomorphic change of coordinates $w_1=\sqrt{\lambda_1} z_1 + i \sqrt{\lambda_2} z_2, w_2=\sqrt{\lambda_1} z_1 -i \sqrt{\lambda_2} z_2,$
i.e., $z_1=\frac{w_1+w_2}{2 \sqrt{\lambda_1}}, z_2=\frac{w_1-w_2}{2 \sqrt{\lambda_2} i}.$ Then $\Ol{E}_t$ is defined as follows in the new coordinates:
$$w_1 w_2=t,~~~~|w_1|^2 +|w_2|^2 \leq \eta^2.$$
This is equivalent $w_2=\frac{t}{w_1}, |w_1|^2 + |\frac{t}{w_1}|^2 \leq \eta^2.$
This yields $\Ol{E}_t$ is biholomorphic to an annulus  $e^{a_1} \leq |w_1| \leq  e^{a_2}$ for some $a_2 > a_1$ depending on $\eta$ and $t$.
Furthermore, by using the $w-$coordinates, we note $\Ol{E}_t$ is covered by a closed strip $\mathcal{S}_t:=\{ \xi \in \mathbb{C}: a_1 \leq \mathrm{Re}(\xi) \leq a_2 \}.$
Indeed the map $\pi(\xi)=(e^{\xi}, te^{-\xi})$ gives a covering map from $\mathcal{S}_t$ to $\Ol{E}_t.$ In particular, writing $\partial E_t$ for the boundary of $\Ol{E}_t,$  $\pi(\{ \mathrm{Re}(\xi) = a_1 \})$ corresponds to one component of $\partial E_t$, and $\pi(\{ \mathrm{Re}(\xi) = a_2 \})$ corresponds to the other.

Note by (\ref{eqnrhoest}), for any point $z =(z_1, z_2) \in \Ol{E}_t,$ we have
$\rho(z) \geq (1-\delta) ||z||^2 +2t +b > b.$
Hence we have $\Ol{E}_t \subset \Ol{\Omega}_{b, 0}$.
On the other hand,  if $z=(z_1, z_2)$ is on the boundary $\partial E_t$ of $\Ol{E}_t,$ it follows from (\ref{eqnrhoest}) that
$$\rho(z) \geq (1-\delta)||z||^2 + 2t +b > (1-\delta)\frac{\eta^2}{ 2 \lambda_1 + 2 \lambda_2} +b. $$
As before, this implies there exists a positive number $A_1$ independent of $t$ such that
\begin{equation}\label{eqnradius}
d(\partial{E_t}, \partial (\Ol{\Omega}_{b, 0} \cap U)) \geq
A_1~\text{for all}~t.
\end{equation}
Note there also exists a positive number $A_2$ independent of $t$  such that
\begin{equation}\label{eqnradius2}
d(E_t, \partial U) \geq A_2~\text{for all}~t.
\end{equation}
Set $A=\mathrm{min}\{ A_1, A_2 \}.$
Let $\epsilon$ be a small positive number such that $V_{\epsilon} \subset U$ and $d(z, 0) < A$ whenever $z \in V_{\epsilon}.$
Write for $0 < t < \widetilde{\epsilon},$ $z_t=(\sqrt{\frac{t}{\lambda_1}}, 0)$ and note $z_t \in \Ol{E}_t$ and $z_t \rightarrow 0$ as $t \rightarrow 0$. Now fix any branch $[g]_q$ with $q \in V_{\epsilon} \cap \Ol{\Omega}_{b, 0}$ as in the assumption of Lemma \ref{lemmaqvepsi}. As $V_{\epsilon} \cap \Ol{\Omega}_{b, 0}$ is connected, we can first extend $[g]_q$ holomorphically along certain curve in $V_{\epsilon} \cap \Ol{\Omega}_{b, 0}$ to  obtain a new branch $[g_0]_{z_t}$  for sufficiently small $0 <t <\widetilde{\epsilon}$.  Then note $\Ol{E}_t$ and thus a small neighborhood of it are connected(but not simply connected).  We can thus extend $[g_0]_{z_t}$ to a multi-valued function, still denoted by $g,$ in a small neighborhood $V$ of $\Ol{E}_t.$ By the boundedness assumption on the holomorphic continuation in $\Ol{\Omega}_{b, 0},$ the norm of any branch of $g$ is bounded by some constant $C>0$. Since $\Ol{E}_t \subset \Ol{\Omega}_{b, 0} \cap U,$ there exists a constant $r_t >0$ depending on $t$ such that the polydisk $\mathbb{P}(z, r_t) \subset \Ol{\Omega}_{b, 0} \cap U$ for any point $z \in \Ol{E}_t$. Then any branch $[g]_z$ of $g$ at $z$ extends to a single-valued holomorphic function, still denoted by $g$, in $\mathbb{P}(z, r_t).$ We apply Cauchy's estimate to $g$ on $\mathbb{P}(z, r_t)$ to obtain for any multiindex $\alpha=(\alpha_1, \alpha_2),$
\begin{equation} \label{eqngald}
|D^{\alpha} g(z)| \leq \frac{C \alpha!}{(r_t)^{|\alpha|} }.
\end{equation}
Let $p$ be a point on $\partial \Ol{E}_t.$ By (\ref{eqnradius}) and (\ref{eqnradius2}), we see that $\mathbb{P}(p, A) \subset \Ol{\Omega}_{b, 0} \cap U.$
Then any branch $[g]_p$ of $g$ at $p$ extends to a single-valued holomorphic function $\hat{g}$ in $\mathbb{P}(p, A).$ As $\Ol{E}_t \cup \mathbb{P}(p, A) \subset \Ol{\Omega}_{b, 0},$ by the boundedness assumption on holomorphic continuation in $\Ol{\Omega}_{b, 0},$ we have $|\hat{g}| \leq C$ in $\mathbb{P}(p, A).$ It then follows from Cauchy's estimate that for any multiindex $\alpha=(\alpha_1, \alpha_2),$
\begin{equation}\label{eqnhatg}
|D^{\alpha} \hat{g}(p)| \leq \frac{C \alpha!}{ A^{|\alpha|}}.
\end{equation}

By the monodromy theorem, the multi-valued $g$ (respectively its derivatives $D^{\alpha} g$) on $\Ol{E}_t$ lifts to a single-valued holomorphic function $h$ (respectively its derivatives $D^{\alpha} h$) in the covering $\mathcal{S}_t$ which is holomorphic in the interior of $\mathcal{S}_t$ and continuous to the boundary $\partial \mathcal{S}_t.$ Moreover, it follows from (\ref{eqngald})  that
$|D^{\alpha} h(\xi)| \leq \frac{C \alpha!}{(r_t)^{|\alpha|}}$ for all $\xi \in \mathcal{S}_t$. Thus $D^{\alpha} h$ is bounded on $\mathcal{S}_t.$ It follows from (\ref{eqnhatg}) that $|D^{\alpha} h(\xi_0)| \leq \frac{C \alpha!}{ A^{|\alpha|}}$ for $\xi_0 \in \partial \mathcal{S}_t.$ We apply the Phragm\'en-Lindel\"of principle to $D^{\alpha} h$ on $\mathcal{S}_t$ to obtain that for all $\xi \in \mathcal{S}_t,$
$|D^{\alpha} h(\xi)| \leq \frac{C \alpha!}{ A^{|\alpha|}}.$
This implies that for every $z \in \Ol{E}_t$, we have
\begin{equation}
|D^{\alpha} g(z)| \leq \frac{C \alpha!}{ A^{|\alpha|}}~\text{for any branch of}~g~\text{at}~z.
\end{equation}
In particular it holds at $z=z_t \in \Ol{E}_t$ and for the branch $[g_0]_{z_t}.$
Hence the Taylor expansion of $[g_0]_{z_t}$ about $z_t$ is convergent in $\mathbb{P}(z_t, A)$, and thus $[g_0]_{z_t}$ extends to a holomorphic function in $\mathbb{P}(z_t, A)$. In particular, $[g_0]_{z_t}$ extends to a holomorphic function $h$ in $V_{\epsilon}$ as $V_{\epsilon} \subset \mathbb{P}(z_t, A)$ for small $t >0.$

\bigskip

Finally, we prove the $\psi-$estimate on $h$. We will also need to shrink $V_{\epsilon}.$ For fixed small $\lambda, \mu> 0,$ we define a complex parameter family (parametrized by $\tau \in \mathbb{C}$) of Riemann surfaces with boundary $\Ol{E}_{\tau},  |\tau| < \mu,$ given by
$$\Ol{E}_{\tau}:=\{z=(z_1, z_2) \in U :\lambda_1 z_1^2 +\lambda_2 z_2^2=\tau;~~~~2 \lambda_1 |z_1|^2 + 2 \lambda_2 |z_2|^2 \leq \lambda^2 \}.$$
Then we apply the same argument as in {\bf (b)} to get the desired region $V_{\epsilon_3}$ and the $\psi-$estimate. And the boundedness follows similarly.

\bigskip

This finishes the proof of Lemma \ref{lemmaqvepsi}.
\end{proof}

Now we let $\{p_1, \cdots, p_k \}$ be the critical points of $\rho$
on $M_b.$ By Lemma \ref{lemmaqvepsi}, we can choose for $1 \leq j
\leq k$ small simply connected neighborhoods $V_{1, j} \ssubset
V_{2, j} \ssubset V_{3, j}$ of $p_j$ such that any branch $[g]_q$ of
$(f, D)$ with $q \in  V_{3, j} \cap \Ol{\Omega}_{b, 0}$ (obtained by
holomorphic continuation along some curve in $\Ol{\Omega}_{b, 0}$)
extends to a holomorphic function in $V_{3, j}$ with
$\psi-$estimate. Moreover, by making $V_{1, j}$ sufficiently small,
we can assume $V_{3,j} \setminus V_{1, j}$ is connected for all $j$.
Let $\delta_1$ be such that if $z_1 \in V_{2, j}$ and $\hat{d}(z_1,
z_2) < 2 \delta_1,$ then $z_2 \in V_{3, j}.$ Now choose $0 <
\epsilon' \ll 1$ such that there is a continuous retract $J$ of
$\left (\Ol{\Omega}_{b-\epsilon', 0} \setminus \cup_{j=1}^k  V_{1,
j}\right)\cup \Ol{\Omega}_{b+\epsilon', 0}$ into
 $\Ol{\Omega}_{b+\epsilon', 0}$
  (See $\S
 3.1$). Here $J$  maps a point in $\Ol{\Omega}_{b+\epsilon', 0}$
 to itself and maps the other points along the orbit of $X_\rho$
 to $M_{b+\epsilon'}$.
 Then by applying a similar argument as in $\S 3.1$, we obtain a
sufficiently small $\epsilon' >0$ such that $(f, D)$ admits
unrestricted holomorphic extension in
 $\Ol{\Omega}_{b-\epsilon', 0} \setminus \cup_{j=1}^k V_{1, j}$ with $\psi-$estimate.
  More precisely, write
 $\mathcal{D}$ for the union of all small holomorphic disks attached
to $M_{b+\epsilon'}$  that can be continuously  deformed to a point
in $M_{b+\epsilon'}$. Choosing $\epsilon'$ small, we make $\Ol{\Omega}_{b-\epsilon',
b+\epsilon'} \setminus \cup_{j=1}^k V_{1, j}$ be
 contained in $\mathcal{D}.$
 We can further make $\epsilon'$ small such that
  for any curve $\sigma$ in $  \Ol{\Omega}_{b-\epsilon', 0}
\setminus \cup_{j=1}^k V_{1, j}$ with $\sigma(0) \in D,$ we can
deform $\sigma$ to a $\widetilde{\sigma}$ in $\Ol{\Omega}_{b}$  through $J$  such that
$\sigma(0)=\widetilde{\sigma}(0)$ and $\sigma(t),
\widetilde{\sigma}(t)$ are sufficiently close for each $t\in [0,1]$.
Moreover, the holomorphic continuation of $(f, D)$ along the two
curves are induced by the same branch at each $t$.

The above means $\Ol{\Omega}_{b-\epsilon', 0} \setminus \cup_{j=1}^k
V_{1, j}$ has the extendable property. Now we make the following
claim:

\medskip
{\bf Claim:} $\Ol{\Omega}_{b-\epsilon', 0}$ has the extendable property.

\medskip
{\bf Proof of Claim:} Let $\gamma$ be any curve in $\Ol{\Omega}_{b-\epsilon', 0}$ with $\gamma(0)
 \in D.$ We first find a $\delta^*$ such that whenever $\gamma(t) \in \Ol{V_{1, j}}$ for some $1
 \leq j \leq k$ and $|t-t'| < \delta^*$ with $t, t' \in [0, 1],$ we have $\gamma(t') \in V_{2, j}.$
 Note if $\gamma([0, 1]) \subset \Ol{\Omega}_{b-\epsilon', 0} \setminus \cup_{j=1}^k V_{1, j}$,
 then the proof is done. Now assume there is some $t_1^*$ such that $\gamma([0, t_1^*)) \subset \Ol{\Omega}_{b-\epsilon', 0} \setminus \cup_{j=1}^k V_{1, j}$ and $\gamma(t_1^*) \in \Ol{V_{1, j_0}}$ for some $1 \leq j_0 \leq k$.  Fix a number $\hat{t}_1$ with $\hat{t}_1 < t_1^*$ and close to $t_1^*$ such that $\gamma(t) \in V_{2, j_0}$ when $\hat{t}_1 \leq t  \leq t_1^*.$
 By the above argument, deforming $\gamma([0, \hat{t}_1])$ if necessary, we can assume there is some $0< t_0 <\hat{t}_1$ such that  $\gamma([0, t_0]) \subset \Ol{\Omega}_{b, 0}$ and $\hat{d}(\gamma(t), \gamma(\hat{t}_1)) < 2 \delta_1$ if $t_0 \leq t \leq \hat{t}_1.$ But since $\gamma(\hat{t}_1) \in V_{2, j_0},$ we have $\gamma(t) \in V_{3, j_0}$ when $t_0 \leq t \leq \hat{t}_1$  (and thus $\gamma(t) \in V_{3, j_0}$ for $t_0 \leq t \leq t_1^*$).
Writing $[f]_{\gamma(t)}$ for the holomorphic continuation of $(f, D)$ along $\gamma,$ we conclude $[f]_{\gamma(t_0)}$
is a branch of $(f, D)$ obtained by holomorphic continuation along a curve in $\Ol{\Omega}_{b, 0}$ to a point in $\Ol{\Omega}_{b, 0} \cap V_{3, j_0}$. It follows from our assumption on $V_{3, j}$ that $[f]_{\gamma(t_0)}$ extends to a holomorphic function $h$ in $V_{3, j_0}$ with $\psi-$estimate and thus $[f]_{\gamma(t)}$ is the restriction of $h$ at $\gamma(t)$ for $t_0 \leq t \leq t_1^*.$ If $\gamma([t_1^*, 1]) \subset V_{3, j_0}$, then clearly $(f, D)$ admits holomorphic continuation along $\gamma$ with $\psi-$estimate. If $\gamma([t_1^*, 1]) \not \subset V_{3, j_0},$ then there exists some $t_2$ with $t_1^* < t_2 < 1$ such that $\gamma([t_1^*, t_2)) \subset V_{2, j_0}$ and $\gamma(t_2)  \in V_{3, j_0} \setminus V_{2, j_0}.$ Since $\gamma(t_1^*) \in \Ol{V_{1, j_0}},$ we must have $|t_2-t_1^*| \geq \delta^*.$ Moreover, by the proceeding argument, we can deform $\gamma([t_0, t_2])$ in $V_{3, j_0}$ with endpoints fixed such that $\gamma$ avoids $V_{1, j_0}$ and we still get the same branch at $\gamma(t_2).$ In other words, we can obtain the same branch of function at $\gamma(t_2)$ by extends $(f, D)$ along a curve in $\Ol{\Omega}_{b-\epsilon', 0} \setminus \cup_{j=1}^k V_{1, j}.$ Next we consider $\gamma([t_2, 1])$. If $\gamma([t_2, 1]) \subset \Ol{\Omega}_{b-\epsilon', 0} \setminus \cup_{j=1}^k V_{1, j},$ then the proof is done again. Otherwise we repeat the above argument for at most $[\frac{1}{\delta^*}]+1$ times to arrive at $\gamma(1)$ and complete the proof of the claim. \endpf

\bigskip

The above claim gives a contradiction to the definition of $b$ and thus {\bf Case III} cannot happen. And it finishes the proof of part (1) and (2) in Theorem \ref{thm33}.

To prove the last part of Theorem \ref{thm33}, we assume there is a holomorphic branch $(h, U)$ of $(f, D)$ in $\mathrm{Reg}(\Ol{\Omega})$ such that $\psi(h(p))=0$ for some $p \in U \setminus M.$ But by part (1) and (2) of Theorem \ref{thm33}, $\psi (h) \leq 0$ on $U$. Since $\psi \circ h$ is subharmonic,  it follows from the maximum principle that $\psi(h) \equiv 0$ in $U$. Let $\sigma$ be the curve in $\mathrm{Reg}(\Ol{\Omega})$ along which we obtain $(h, U)$ by applying holomorphic continuation to $(f, D).$  And we write $\{(h_l, U_l) \}_{l=0}^{\nu}$ for the holomorphic continuation of $(f, D)$ along $\sigma$ with $U_{\nu}=U$ and $h_{\nu}=h.$ Since $h_l=h_{l+1}$ on $U_l=U_{l+1},$ we have each $\psi (h_l)$ attains its maximum at an interior point and thus $\psi (h_l) \equiv 0$ on $U_l$ for all $l$. This in particular holds for $l=0$ and note $U_{0} \cap \Ol{\Omega}_{-\epsilon_2, 0}$ has a component containing $\sigma(0) \in D$. Then by Corollary \ref{coromegae}, we have $\psi(g) \equiv 0$ for every holomorphic branch $(g, V)$ of $(f, D)$ in $\Ol{\Omega}_{-\epsilon_2, 0}$. This in particular implies $\psi(f) \equiv 0$ on $D$. Moreover, as every curve in $\mathrm{Reg}(\Ol{\Omega})$ from $(f, D)$ must pass through $\Ol{\Omega}_{-\epsilon_2, 0}$, it also follows that $\psi(g) \equiv 0$ for every holomorphic branch $(g, V)$ of $(f, D)$ in $\mathrm{Reg}(\Ol{\Omega})$. This establishes Theorem \ref{thm33}.







\subsection{Proof of Theorem \ref{theorem2}}

In this section, we use Theorem \ref{thm33-preliminary}(or Theorem \ref{thm33}) to prove Theorem \ref{theorem2}. Write $n$ for the complex dimension of $\Omega$ and $\mathbb{B}^n$ for the complex unit ball in $\mathbb{C}^n.$ By assumption, for any $p_0 \in M$, there is a small open piece $D$ containing $p_0$ of $M$ and a smooth CR diffeomorphism $f$ from $D$ to $\partial \mathbb{B}^n.$ Fix such a CR function element $(f, D).$
By Alexander \cite{Al}, any two CR diffeomorphisms from a connected open
piece of $M$ to $\partial \mathbb{B}^n$ must differ by an
automorphism of $\mathbb{B}^n$ for $n\ge 2$. This implies $(f, D)$
admits CR continuation $[f]_{\sigma(t)}$ along any path $\sigma$ in
$M$ with the image of $[f]_{\sigma(t)}$ in $\partial \mathbb{B}^n$
by Burns and Shnider \cite{BSc}.  In particular, taking
$\psi(z)=||z||^2 -1$ for $z \in \mathbb{C}^n,$ $(f, D)$ admits
unrestricted (bounded) CR continuation in $M$ with $\psi-$estimate.
Furthermore, for any $p \in M$ and any two CR branches $[f_1]_p$ and
$[f_2]_p$ of $(f, D)$, we must have $f_2=G \circ f_1$ near $p$  for
some automorphism $G$ of $\mathbb{B}^n.$ Let $\Gamma_p$ be the
collection of all such $G's$. Apparently, $\Gamma_p $ is a subgroup
of the automorphism group of the unit ball.  By uniqueness of CR
functions extending to the pseudoconvex side, we have
$\Gamma:=\Gamma_p$ is independent of $p\in M$.

By Theorem \ref{thm33}, $(f, D)$ admits unrestricted holomorphic continuation in $\mathrm{Reg}(\Ol{\Omega})$ with $\psi-$estimate. Note our choice of $\psi$ is strongly pseudoconvex in $\mathbb{C}^n$ and $f$ is not constant. By Remark \ref{remarkstrict}, the image of any holomorphic branch $(g, V)$  must stay in $\mathbb{B}^n$ if $V \subset \mathrm{Reg}(\Omega).$ Let $\Ol{\Omega}_{-\epsilon_2, 0}$ be as in the
proof of Theorem \ref{thm33}. It follows readily from the construction of the holomorphic continuation
 in $\Ol{\Omega}_{-\epsilon_2, 0}$ that any holomorphic branches $[g_1]_q$ and
$[g_2]_q$ of $(f, D)$ at some $q \in \Ol{\Omega}_{-\epsilon_2, 0}$
must satisfy $g_2=G \circ g_1$  for some automorphism $G\in \Gamma$.
We further claim
the following statement.

\medskip
{\bf Claim:} For any two branches $[h_1]_q$ and $[h_2]_q$ of $(f, D)$ at some $q \in \mathrm{Reg}(\Omega),$
 it must hold that $h_2=G \circ h_1$  for some automorphism $G\in
 \Gamma$ with $G$ depending on $h_1$ and $h_2$.

\medskip
{\bf Proof of Claim:} Let $\gamma_1$ and $\gamma_2$ be the two curves with $p_0=\gamma_1(0)= \gamma_2(0)
 \in D$ along which we obtain the two branches $[h_1]_q$ and $[h_2]_q$ respectively  by holomorphic
  continuation of $(f, D).$  Set  the path $\alpha$ to be the sum of $\gamma_1$ and
  the reverse $-\gamma_2$ of $\gamma_2: \alpha= \gamma_1 + (-\gamma_2).$ More precisely,
   $\alpha(t)=\gamma_1(2t)$ if $0 \leq t \leq \frac{1}{2}$ and $\alpha(t)=\gamma_2(2-2t)$
   if $\frac{1}{2} \leq t \leq 1.$ We pause to recall the following proposition  proved by Huang-Ji in
\cite{HJ}, (see Lemma 5.2 \cite{HJ}), which in particular implies
that the branches obtained by continuing along loops  in
$\overline{\Omega}_{r_2,0}$ based on $p\in
\overline{\Omega}_{r_1,0}$ with $r_1>r_2$ are exactly those obtained
long loops inside $\overline{\Omega}_{r_1,0}$ based on $p$:
\begin{proposition}\label{corogammatil}
Let $\gamma$ be a curve in $\mathrm{Reg}(\Ol{\Omega})$ with
$\gamma(0) \in D$ and $\gamma(1) \in \Ol{\Omega}_{-\epsilon_2, 0}.$
Then $\gamma$ can be continuously deformed  inside
$Reg(\overline\Omega)$ with endpoints fixed to a new curve
$\widetilde{\gamma} \subset \Ol{\Omega}_{-\epsilon_2, 0}$ and thus
we obtain the same branch at $\gamma(1)=\widetilde{\gamma}(1)$ by
continuing  $(f, D)$ along the these two curves.
\end{proposition}


By Proposition \ref{corogammatil},
    $\alpha$ can be deformed with endpoints fixed to a curve in $\Ol{\Omega}_{-\epsilon_2, 0}$
    without changing the holomorphic branch we obtain at the endpoint $p_0$. Then by the preceding discussion,
    writing $(\hat{f}, U)$ with $p_0 \in U$ for the holomorphic branch we get by
     holomorphic continuation of $(f, D)$ along $\alpha$, we must have $\hat{f}=H \circ f$ on $V$ for
     some automorphism $H\in\Gamma$. This yields that if we apply holomorphic continuation
     to $(H \circ f, D)$ along $\gamma_2,$ we get the branch $[h_1]_q$ at $q$. Consequently,
     we have $h_2=G \circ h_1$ near $q$ where $G=H^{-1}.$ This proves the claim. \endpf

\medskip
 Now we define the complex analytic hyper-variety
$E\subset \mathrm{Reg}(\Omega)$ to be  such that,  for any branch
$(f^*,U^*)$ of $(f,D)$, $E\cap U^*$ is  the
zero of the Jacobian of $f^*$. Then we see from the above claim that $E$ is well-defined
and is independent of the choice of the chosen branch. Since $E\cap
\Omega_{-\epsilon_2,0}=\emptyset$, we see that $E=\emptyset$.
Namely, $f^*$ is always a local biholomorphism.

We now define the hyperbolic metric $\omega_0$ on
$\mathrm{Reg}(\Omega)$ in the following way. Writing
$\omega_{\mathbb{B}^n}$ for the Bergman (hyperbolic) metric on
$\mathbb{B}^n,$ for any holomorphic branch $(g, V)$ of $(f, D)$ in
$\mathrm{Reg}(\Omega),$
 we define $\omega_0=g^*(\omega_{\mathbb{B}^n})$ on $V$. The above claim guarantees that  $\omega_0$
 is a K\"ahler metric which is  independent of the choice of $(g, V)$ as the hyperbolic metric
  on $\mathbb{B}^n$ is invariant under
 automorphisms. And thus the metric $\omega_0$ is well-defined on $\mathrm{Reg}(\Omega)$ . Finally
 we notice that
  for any $p \in M,$ there are a neighborhood $W_p$ of $p$ in $\mathrm{Reg}(\Ol{\Omega})$ and a
  smooth
   diffeomorphism $F$ from  $W_p$
    to a certain open subset $W'_p$
   of  $\Ol{\mathbb{B}^n}$ such that (i) $\partial_0 W_p:=W_p\cap M$ is an open
   subset of $M$ containing $p$ and $\partial_0 W'_p:=W'_p\cap
   \partial{\mathbb{B}^n}$ is an open subset of the unit sphere
   $\partial{\mathbb{B}^n}$ (ii) $F$ is  CR diffeomorphism from $\partial_0 W_p$
    to $\partial_0 W'_p$ and extends
   to a biholomorphism from $W_p\cap \Omega$ to $W'_p\cap {\mathbb
   B}^n$ (iii)
    $\omega_0=F^*(\omega_{\mathbb{B}^n})$
on
    $W_p \setminus M$.
This implies $\omega$ is complete near $M$. This finishes the proof
of Theorem \ref{theorem2}. $\endpf$

\bigskip

To end this section, we give a result on the uniqueness of the
hyperbolic metric on $\mathrm{Reg}(\Omega).$ Fix $p_0 \in M.$ We
emphasize that an important feature of the metric $\omega_0$ is that
there is a small neighborhood $W$ of $p_0$ in
$\mathrm{Reg}(\Ol{\Omega})$ such that
$\omega_0=F^*(\omega_{\mathbb{B}^n})$ on $W \setminus M$, where $F$
is a CR diffeomorphism from $W\cap M$ to an open piece of $\partial
{\mathbb B}^n$ that extends to a holomorphic embedding from
$W\setminus M$ into  ${\mathbb{B}^n}$. This property makes the
metric unique in the case of complex dimension at least two. Indeed,
we will establish a uniqueness result for the hyperbolic metric
under a weaker boundary condition.
Let $\hat{\Omega}$ be a connected open subset of
$\mathrm{Reg}(\Omega).$ Assume $\hat{\Omega}$ has $\hat{M}\subset M$
as part of its smooth boundary with  $p_0 \in \hat{M}.$


\begin{definition}\label{defncomplete}
Let $\omega$ be a metric on $\hat{\Omega}.$ We say $\omega$ is boundary complete at $p_0 \in M$ if
there exists a small simply connected open subset $U$ of
$\hat{\Omega}$ with $p_0$ on its smooth part of boundary, a point
$q_0 \in U$ and an open piece $\Gamma \subset \mathbb{S}^{2n-1}
\subset T_{q_0} U$ such that for every $v \in \Gamma$, there is a
geodesic (with respect to $\omega$) $\gamma_v : [0, \infty)
\rightarrow U$ with $\gamma_v(0)=q_0, \gamma'_v(0)=v$ and
$\gamma_v(t)$ converges to some point on $M$ as $t$ goes to
infinity. Moreover, writing $V$ for  the union of the geodesics $\{ \gamma_v(t) \}_{v \in \Gamma},$ there is a neighborhood $V_0$ of $p_0$ in $\mathrm{Reg}(\Ol{\Omega})$ such that $V_0 \setminus M \subset V.$
\end{definition}

We note the metric $\omega_0$ we construct is boundary complete in the sense of Definition \ref{defncomplete}.
We will show that a hyperbolic metric must be unique if it is boundary complete in the sense of Definition \ref{defncomplete} and the dimension is at least two.

\begin{proposition}\label{propcomplete}
Let $\Omega$ be the Stein space as in Theorem \ref{theorem2} of dimension at least two,  and $\hat{\Omega}, \hat{M}, $and $p_0$ be as above. There exists a unique metric on $\hat{\Omega}$ up to scaling with the following property.
\begin{enumerate}[(a)]
\item $\omega$ is hyperbolic on $\mathrm{Reg}(\hat\Omega).$
\item $\omega$ is boundary complete at $p_0.$
\end{enumerate}
\end{proposition}

\begin{proof}
The existence of the desired metric follows from Theorem \ref{theorem2}.  Indeed, we can obtain the
 desired metric $\omega$ by restricting the metric $\omega_0$ (we constructed in Theorem \ref{theorem2})
 to $\hat{\Omega}.$

We now prove the uniqueness part. By assumption, we can assume the statement
 in Definition \ref{defncomplete} holds and let $U, q_0, \gamma_{v}, V$ be as there.
 First since $\omega$ is hyperbolic, there exists a small open subset $U_0$ of $U$ and a holomorphic
  isometric map $F: (U_0, \lambda \omega) \rightarrow (\mathbb{B}^n, \omega_{\mathbb{B}^n})$ for some
   normalizing constant $\lambda > 0.$ By a classical theorem of Calabi \cite{Ca}, $F$ extends holomorphically
   along any path in $\hat{\Omega}.$ By the monodromy theorem, it extends to a holomorphic isometric
    immersion, still denoted by $F$, from $(U, \lambda \omega)$ to $(\mathbb{B}^n, \omega_{\mathbb{B}^n}).$
    In particular, we have $F^*(\omega_{\mathbb{B}^n})= \lambda \omega$ on $U$.

Recall $V$ denotes the union of traces of all $\gamma_v, v \in \Gamma.$ It is a consequence of the classical
  Rauch comparison that manifolds of nonpositive curvature have no conjugate points and the exponential map
   is a local diffeomorphism. We thus have $V$ is an open subset of $U$. Moreover,  $V$ has a piece of smooth boundary $M_0 \subset M.$ By composing with an automorphism
    of $\mathbb{B}^n,$ we assume $F(q_0)=0 \in \mathbb{B}^n.$

\medskip
{\bf Claim:} $F$ maps every geodesic $\gamma_v$ to a (straight) ray with initial point $0$ in $\mathbb{B}^n.$ Moreover, $F$ is an embedding from $V$ to $\mathbb{B}^n,$ and the image of $V$ under $F$ is a conic region $W$ in $\mathbb{B}^n.$

\medskip
{\bf Proof of Claim:} First since $F$ is local isometric and $F(q_0)=0$, $F(\gamma_v)$ must be a geodesic in $(\mathbb{B}^n, \omega)$ with intial point at $0$, which is a ray in this case. Furthermore, $F$ maps two distinct geodesics $\gamma_{v_1}$ and $\gamma_{v_2}$ to two different rays in $\mathbb{B}^n$.
Also we claim $F$ maps two different points $\gamma_v(t_1)$ and $\gamma_v(t_2),  0 \leq t_1 < t_2,$ on a geodesic to two different points on a ray. Otherwise, $F \circ \gamma_v (t_1)=F \circ \gamma_v (t_2)$ and there must exist  some $t_1 < t_3 < t_2$ such that $F \circ \gamma_v (t)$ is not embedding near $t_3$, a plain contradiction to the fact that $F$ is a local embedding.
This implies $F$ is an emdedding on $V$. Recall $V$ is the union of the traces of all geodesics. The other part of the claim follows as well. \endpf

\medskip
Now by the above claim, $F$ is a biholomorphic map from $V$ to $W$, and maps a geodesic in $V$ (which converges to a point on $M$) to a ray in $W$ which converges to a point on $\partial \mathbb{B}^n,$ respectively.
By a theorem of Forstneri\v{c}-Rosay  (Page 239 in \cite{FR}), as $M$ and $\partial \mathbb{B}^n$ are strongly pseudoconvex, $F$ extends continuously up to a piece of smooth boundary $N \subset M$ of $V$ such that $F(N) \subset \partial \mathbb{B}^n.$ Since $M$ is spherical, it then follows from Alexander \cite{Al} that $F$ is a smooth CR diffeomorphism from $N$ to $F(N).$ We run our process to $(F, N)$ in Theorem \ref{theorem2} to obtain the metric $\omega_0$ on $\mathrm{Reg}(\Omega).$ We claim this metric $\omega$ is identical with $\omega_0$ on $\hat{\Omega}$ after an appropriate scaling. Indeed, recall to obtain $\omega_0,$ we extend to $(F, N)$ holomorphically the pseudoconvex side and apply holomorphic continuation along any path to get a new branch $(g, V_1).$ Then pull back the hyperbolic metric $\omega_{\mathbb{B}^n}$ on $\mathbb{B}^n$ by $g$ to get the metric $\omega_0$ on $V_1$. But the holomorphic extension of $(F, N)$ to the pseudoconvex side is precisely $(F, V)$ and thus the metric $\omega_0$ on $V$ is given by $\omega_0=F^*(\omega_{\mathbb{B}^n}).$ Thus $\omega_0=\lambda \omega$ on $V$. Since both $\lambda \omega$ and $\omega_0$ are real analytic and $\hat{\Omega}$ is connected, we finally conclude $\lambda \omega$ is the restriction of $\omega_0$ on $\hat{\Omega}.$ This proves the uniqueness of the metric.
\end{proof}

\medskip

We revisit Example \ref{epdisk} to show the uniqueness in Proposition \ref{propcomplete} fails in one dimensional case.

\begin{examplenonum} (Example \ref{epdisk} in $\S1$)
Let $\Delta$ be the unit disk and $X$ be the Riemann surface with singularity in $\mathbb{C}^2$ give by $\{ (z, w) \in \mathbb{C}^2: w^2=z^3\} \cap \Delta^2.$ It is indeed the image of the map $t \rightarrow (t^2, t^3), t \in \Delta$ and has an isolated singularity $(0, 0).$
Note $X^*:=X \setminus \{ (0, 0) \}$ is biholomorphic to punctured disk $\Delta^*:=\Delta \setminus \{ 0 \}.$ Then the hyperbolic metric on $\Delta^*$ induces a
metric $\omega_1$ on $X^*$ and the Bergman metric on $\Delta^*$ induces a metric $\omega_2$ on $X^*$. Note $\omega_1$ is complete near the singularity  and $\omega_2$ is not. We claim that both $\omega_1$ and $\omega_2$ satisfy the boundary completeness in Definition \ref{defncomplete}. Indeed, $\omega_2$ is identical with the Bergman metric of $\Delta$ near boundary and thus satisfies the boundary completeness in Definition \ref{defncomplete}. To understand $\omega_1$, we look at the covering map $\pi (\xi)=e^{i \xi}$ from the upper half plane $\mathcal{H}=\{ \xi \in \mathbb{C}: \mathrm{Im}~\xi > 0 \}$ to $\Delta^*.$ Note $\pi$ maps $\{ \xi \in \mathcal{H}: 0 <  \mathrm{Re} ~\xi  < 2 \pi \}$ to an open subset of $\Delta^*$ whose boundary contains an open piece of the circle. The implies the induced hyperbolic metric on $\Delta^*$ and thus $\omega_1$ on $X^*$ satisfies the boundary completeness at some boundary point in the sense Definition \ref{defncomplete}. Hence we have two distinct metrics on $X^*$ that satisfy both (a) and (b).
\end{examplenonum}

\section{Proof of Theorem \ref{theorem3}}

The section is devoted to establishing Theorem \ref{theorem3}. Let $\Omega$ be a
 Stein space with possibly isolated singularities and write $\mathrm{Reg}(\Omega)$ for its regular part.
 Write $\Lambda^n(\mathrm{Reg}(\Omega))$ for the space of the holomorphic $n-$forms
  on $\mathrm{Reg}(\Omega)$ and define the Bergman space
  of $\Omega$ to be $A^2(\Omega)=\{\alpha \in \Lambda^n(\mathrm{Reg}(\Omega)):
  |\int_{\mathrm{Reg}(\Omega)} \alpha \wedge \Ol{\alpha}| < \infty \}.$ Assume $A^2(\Omega)$ is not empty. Then $A^2(\Omega)$ is a Hilbert space with inner product:
$$\langle \alpha, \beta \rangle :=\frac{1}{i^{n^2}} \int_{\mathrm{Reg}(\Omega)}
 \alpha \wedge \Ol{\beta}~\text{for}~\alpha, \beta \in \Lambda^n(\mathrm{Reg}(\Omega)).$$

Let $\{ \alpha_j \}_{j=1}^{m}$ be an orthonormal basis of $A^2(\Omega), m \leq \infty,$ and define the Bergman kernel to be
$K_{\Omega}=\sum_{j=1}^{m} \alpha_j \wedge \Ol{\alpha}.$ In a local holomorphic coordinate chart $(U, z)$ of $\mathrm{Reg}(\Omega),$ we have
$$K_{\Omega}=k(z, \Ol{z}) dz_1 \wedge \cdots dz_n \wedge d\Ol{z_1} \wedge \cdots \wedge d\Ol{z_n}~\text{in}~U.$$
Here $k(z, \overline{z})=\sum_{j=1}^{m} |a_j|^2$ if $\alpha_j=a_j dz_1 \wedge \cdots \wedge dz_n, j \geq 1,$ in the local coordinate chart.
Assume $K_{\Omega}$ is nowhere zero on $\mathrm{Reg}(\Omega).$ We define a Hermtian $(1, 1)-$form to be $\omega^B=\sqrt{-1} \partial \Ol{\partial} \log k(z, \Ol{z}).$ We call $\omega^B$ the Bergman metric if it induces a (positive definite) metric on $\mathrm{Reg}(\Omega).$

To further study the Bergman metric, we pause to recall the definition of holomorphic maps from a complex manifold $X$ to the infinite dimensional projective space
$\mathbb{P}^{\infty}.$ Let $F$ be a map from $X$ to $\mathbb{P}^{\infty}.$ We say $F$ is holomorphic if for any $p \in X,$ there is a local holomorphic coordinate chart $(U, z)$ with $p \in U$ and a set of holomorphic functions $\{ f_j\}_{j=1}^{\infty}$ on $U$ such that
\begin{enumerate}
\item The set of functions $\{ f_j \}_{j=1}^{\infty}$ is base point free, i.e., they have no common zeros.

\item The infinite sum $\sum_{j=1}^{\infty} |f_j|^2$ converges uniformly on every compact subsets in $U.$

\item $F=[f_1, \cdots, f_n, \cdots ]$ on $U$.
\end{enumerate}
If $X$ is equipped with a K\"ahler metric $\omega,$ we further say
$F$ is isometric if $$\omega=\sqrt{-1}\partial \Ol{\partial} \log
(\sum_{j=1}^\infty |f_j(z)|^2)$$ on $U$. Let $\Omega$ be as above,
equipped with the Bergman metric $\omega^B$ and $\{ \alpha_j
\}_{j=1}^{m}$ an orthonormal basis of $A^2(\Omega).$ Then it
induces a natural holomorphic map $F$ from $\mathrm{Reg}(\Omega)$ to
$\mathbb{P}^{\infty}$ given by
$$F=[\alpha_1, \cdots, \alpha_j, \cdots].$$
Here if $m$ is finite, we then add zero components to make the target be $\mathbb{P}^{\infty}.$ The right hand side of the above equation is understood as follows. In a local holomorphic coordinates chart $(U, z): [\alpha_1, \cdots, \alpha_j, \cdots]=[a_1(z), \cdots, a_j(z),\cdots]$ if $\alpha_j =a_j dz_1 \wedge \cdots \wedge dz_n$. Note this definition is independent of the choices of coordinates. Indeed, if $\alpha_j=b_j dw_1 \wedge \cdots \wedge dw_n$ in another coordinates chart $(V, w),$ then on $U \cap V$ every $b_j$ just differs from $a_j$ by the same factor, i.e., the Jacobian of the change of coordinates. We also remark that $\{ \alpha_j\}_{j=1}^m$ must be base point free by the fact that $K_{\Omega}$ is nowhere zero.  We will call  $F$ the Bergman-Bochner map from $\mathrm{Reg}(\Omega)$ to $\mathbb{P}^{\infty}$ and denote it by $\mathcal{B}_{\Omega}.$ We are now at the position to
prove Theorem \ref{theorem3}.


\medskip
{\bf Proof of Theorem \ref{theorem3}:} The {\it only-if} part is
trivial and we only prove the converse statement. Assume
$(\mathrm{Reg}(\Omega), \omega^B)$ has constant holomorphic
sectional curvature $\lambda.$ {\it A Priori} $\lambda$ can be of
any sign or zero. Write $(X_0, \omega_{st})$ for the space form
where the metric is normalized so that the holomorphic sectional
curvature equals $\lambda.$ More precisely,
\begin{itemize}
\item If $\lambda< 0$, then we let $X_0=\mathbb{B}^n$ and $\omega_{st}$ the (suitably normalized) Poincar\'e metric;

\item If $\lambda=0,$ then we let $X_0=\mathbb{C}^n$ and $\omega_{st}$ the standard Euclidean metric;

\item If $\lambda >0,$ then we let $X_0=\mathbb{P}^n$ and $\omega_{st}$ the (suitably normalized) Fubini-Study metric.
\end{itemize}

Then by assumption, there is a connected open subset $U$ of $\mathrm{Reg}(\Omega)$ and a holomorphic isometric map $f:(U, \omega^B) \rightarrow (X_0, \omega_{st}).$ By a classical theorem of Calabi \cite{Ca}, $f$ extends holomorphically along any path $\gamma \in \mathrm{Reg}(\Omega)$ with $\gamma(0) \in U.$ In other words, $f$ extends to a possibly multi-valued map $F$ from $\mathrm{Reg}(\Omega)$ to $M_0.$ We will prove $F$ must be indeed single-valued.


We will need the following two lemmas. Consider the case $\lambda > 0.$ Then $X_0$ is the projective space $\mathbb{P}^n$ and $\omega_{st}$ is the normalized Fubini-Study metric. That is, writing $[\eta]=[\eta_0, \eta_1, \cdots, \eta_n]$ for the homogeneous coordinates of $\mathbb{P}^n, \omega_{st}$ is given by $\omega_{st}=\mu \sqrt{-1} \partial \overline{\partial} \log ||\eta||^2$ for some constant $\mu > 0.$

\begin{lemma}\label{lemma42}
In the case $\lambda >0$ as above, we must have $\mu$ equals some positive integer $m$.
\end{lemma}

{\bf Proof of Lemma \ref{lemma42}:} By shrinking $U$, we assume $U$ is contained in some holomorphic coordinate chart. And by a holomorphic change of coordinates, we assume $z(p)=0$ for some $p \in U.$
Recall $K(z, \Ol{z})=k(z, \Ol{z})dz_1 \wedge \cdots \wedge  dz_n \wedge d\Ol{z_1} \wedge
\cdots \wedge d\Ol{z_n}$ on $U$.
Without loss of generality, assume $f(0)=[1, 0, \cdots, 0].$ Shrinking $U$ further if necessary,
assume $f(U)$ is contained in the Euclidean cell $V_0=\{[\eta_0, \cdots, \eta_n]: \eta_0 \neq 0\}$ of $\mathbb{P}^n.$ Write
 $g$ for the inverse of $f$ defined from a small neighborhood $V$ of $f(0)$ to $U.$ Note $g$ is also
 isometric: $g^*(\omega^B)=\omega_0$. This yields if we write $\xi=(\xi_1, \cdots, \xi_n)$ for the nonhomogenous coordinates of $\mathbb{P}^n$  on $V_0$, then
$$\mu \partial \Ol{\partial} \log (1+|| \xi || ^2)=\partial \Ol{\partial} \log k(g(\xi), \Ol{g(\xi)})~\text{on}~V.$$
Thus the difference of  $\log (1+|| \xi || ^2)^{\mu}$ and $\log k(g(\xi), \Ol{g(\xi)})$ equals to a pluriharmonic function $\psi$. Shrinking $V$ if necessary, assume $\psi=h+\Ol{h}$ for a holomorphic function $h$ on $V$.  We thus conclude
$$k(g(\xi), \Ol{g(\xi)})|e^{h(\xi)}|^2=(1+|| \xi ||^2)^{\mu}.$$
We then use the definition of $k(z, \Ol{z})$ and apply the Taylor expansion of $(1+x)^{\mu}$ on $\{ x \in \mathbb{R} : |x| < 1 \}$ to get
$$\sum_{j=1}^{\infty} |e^{h(\xi)} (a_j\circ g)(\xi)|^2= 1+\sum_{k=1}^{\infty} \frac{\mu(\mu-1) \cdots (\mu-k+1)}{k!} ||\xi||^{2k} $$
Fix $l \geq 1.$ Apply $\frac{\partial^{2l}}{\partial \xi_1^l \partial \Ol{\xi}_1^l}$ to both sides of the equation and evaluate at $\xi=0$. The left hand side always gives a nonnegative number. But if $\mu$ is non-integer, the right hand side will give a negative number for some sufficiently large $l$, a plain contradiction. Hence $\mu$ must be a positive integer. \endpf


\begin{lemma}\label{lemmaisom}
The normalized space form $(X_0, \omega_{st})$ can be isometrically embedded into $\mathbb{P}^{\infty}$ by a holomorphic map $\mathcal{B}$.
\end{lemma}

{\bf Proof of Lemma \ref{lemmaisom}:}
We start with the case when $\lambda >0.$ By Lemma \ref{lemma42}, we have in this case $X_0=\mathbb{P}^n$ and $\omega_{st}=m \partial \Ol{\partial} \log ||\eta||^2$ for some positive integer $m$. We know $(M_0, \omega_{st})$ can be holomorphically isomerically embedded into  $(\mathbb{P}^N, \omega_{FS})$ by the $m^{\mathrm{th}}$ Veronese embedding for some appropriate $N$. Here $\omega_{FS}$ is the standard Fubini-Study metric on $\mathbb{P}^N.$ Furthermore,
$(\mathbb{P}^N, \omega_{FS})$ can be canonically embedded into $\mathbb{P}^{\infty}$ (by adding zero components). This finishes the proof of Lemma \ref{lemmaisom} in the positive $\lambda$ case.

We next consider the case $\lambda < 0.$ In this case, $X_0=\mathbb{B}^n$ and $\omega_{st}=-\mu \partial \Ol{\partial} \log (1-||w||^2)$ for some $\mu >0,$ where $w=(w_1, \cdots, w_n)$ is the coordinates of $\mathbb{B}^n \subset \mathbb{C}^n.$ We use the Taylor expansion of $(1-x)^{-\mu}$ about the point $x=0$ on $\{ x \in \mathbb{R}: |x| < 1\}$ to obtain
$$(1-||w||^2)^{-\mu}= 1+ \sum_{k=1}^{\infty} \frac{\mu(\mu+1) \cdots (\mu+k-1)}{k!} ||w||^{2k}.$$
The right hand side converges uniformly on compact subsets of $\mathbb{B}^n.$
This implies there exists a sequence of monomials $\{P_j \}_{j=1}^{\infty}$ in $w$ such that $(1-||w||^2)^{-\mu}=\sum_{j=1}^{\infty}|P_j(w)|^2$, which converges uniformly on compact subsets of $\mathbb{B}^n.$ This leads to a natural isometric map $\mathcal{B}=[P_1, \cdots, P_j, \cdots]$ from $(X_0, \omega_{st})$ to $\mathbb{P}^{\infty}.$ It is clear that $F$ is an embedding as the $w_j'$s, $1 \leq j \leq n,$ (with appropriate coefficients) are among the $P_j'$s.

The proof for the case $\lambda=0$ is similar. Note $\omega_{st}=\partial \Ol{\partial} ||w||^2= \partial \Ol{\partial} \log (e^{||w||^2})$ where $w$ is the coordinates of $X_0=\mathbb{C}^n,$ and we use the Taylor expansion of $e^{||w||^2}$ at $w=0$ to obtain an isometric embedding from $X_0$ to $\mathbb{P}^{\infty}.$
This establishes Lemma \ref{lemmaisom}. \endpf

\medskip
Now $U \subset \mathrm{Reg}(\Omega)$ can be isometrically embedded into $\mathbb{P}^{\infty}$ in two ways: (1). $\mathcal{B}_{\Omega}$ embeds $(U, \omega^B)$ isometrically into $\mathbb{P}^{\infty},$ and (2). $(U, \omega^B)$ is isometrically embedded into $\mathbb{P}^{\infty}$ by $\mathcal{B} \circ f.$ By Calabi's theorem \cite{Ca}, there is a rigid motion of $T$ of $\mathbb{P}^{\infty}$ such that
\begin{equation}\label{eqntb}
T \circ \mathcal{B}_{\Omega}= \mathcal{B} \circ f~\text{in}~U.
\end{equation}
As mentined above, $(f, U)$ can be extended holomorphically along any path by \cite{Ca}.
Let $\gamma_1, \gamma_2$ be two curves connecting $p \in U$ to some point $q \in \mathrm{Reg}(\Omega).$ Write $(f_1, V)$ and $(f_2, V)$ for the two holomorphic branches obtained from holomorphic continuaion of $(f, U)$ along the two curves. By the uniqueness of holomorphic continuation, we know (\ref{eqntb}) is preserved along the continuation. Thus we have
$$T \circ \mathcal{B}_{\Omega}= \mathcal{B} \circ f_1= \mathcal{B} \circ f_2~\text{in}~V.$$
But note $\mathcal{B}_{\Omega}, \mathcal{B}$ and $T$ are all embeddings.
This implies $f_1=f_2$ in $V$. Hence $(f, U)$ extends a well-defined
holomorphic map $F$ on $\mathrm{Reg}(\Omega)$, which is local
isometric. By (\ref{eqntb}) again, we have $T \circ
\mathcal{B}_{\Omega}= \mathcal{B} \circ F$ on
$\mathrm{Reg}(\Omega).$ This yields that $F$ is an embedding map.

We next prove that the singular set $\mathrm{sing}(\Omega)$ of $\Omega$ is empty. Suppose not and let $p \in \mathrm{sing}(\Omega).$ First note we can assume $(\Omega, p)$ is embedded into some $\mathbb{C}^K$ with $z(p)=0.$ For a sufficiently small $\epsilon,$ consider the link $N_{\epsilon}=\{z \in \mathbb{C}^K: ||z|| =\epsilon \} \cap \Omega$ at $p$ and $D_{\epsilon}=\{z \in \mathbb{C}^K: ||z|| < \epsilon \} \cap \Omega.$ Note $F$ is a holomorphic map in $D_{\epsilon}-\{p \}$ and extends holomorphically across $N_{\epsilon}.$ Recall $F$ is an embedding. Then $\widetilde{N}=F(N_{\epsilon})$ is a (connected) closed strongly pseudoconvex hypersurace in $X_0$ (See Milnor \cite{Mi}) and $F(D_{\epsilon}-\{p\})$ is contained in some pseudoconvex domain  $\widetilde{D}$ in $X_0$ which has smooth boundary $\widetilde{N}.$ We remark that when $X_0=\mathbb{P}^n, \widetilde{D}$ is defined to be the set of points in $\mathbb{P}^n \setminus \widetilde{N}$ that can be path-connected to a small analytic disk attached to $\widetilde{N}.$ Clearly $\widetilde{D}$ defined in this way is connected and open. To show $\widetilde{D}$ has $\widetilde{N}$ as its smooth boundary, we observe any point $q$ that is on the pseudoconcave side of $\widetilde{N}$ and close to $\widetilde{N}$ cannot be contained in $\widetilde{D}.$ Indeed, suppose $q \in \widetilde{D}.$ Then we find a point
$q'$ on the pseudoconvex side to $\widetilde{N}$ and close to $q.$ Moreover, pick a short smooth curve $\gamma_1$ from $q$ and  $q'$ that
cuts $\widetilde{N}$ transversally. On the other hand, there is a curve $\gamma_2$ in $\widetilde{D}$ from $q'$ to $q$.  Set $\gamma=\gamma_1+ \gamma_2.$ By the construction of $\gamma,$ the intersection number of $\gamma$ and $\widetilde{N}$ is $\pm 1.$ By the simply connectedness of $X_0,$ however, this interesction number must be zero, a plain contradiction. This proves $\widetilde{D}$ is a pseudoconvex domain with  smooth boundary $\widetilde{N}.$

Since $p$ is normal, note we can assume $F$ is the restriction of some holomorphic map $\hat{F}$ on $\{ z \in \mathbb{C}^K: ||z|| < \epsilon\}$ to $D_{\epsilon}.$ This follows from the property of normal singularities if $X_0$ is $\mathbb{B}^n$ or $\mathbb{C}^n,$ and needs some justification when $X_0$ is $\mathbb{P}^n.$ In this case, we use a theorem of Takeuchi \cite{T} to see $\widetilde{D}$ is Stein, which can be embedded as a submanifold of some $\mathbb{C}^{K'}$. Hence $F$ can be regarded
as a map to $\mathbb{C}^{K'},$ and again by the property of normal singularities, $F$ extends continuously across $p$. Furthermore, we can shrink $\epsilon$ to make $\widetilde{D}$ be contained in one Euclidean cell of $\mathbb{P}^n$ and the observation follows readily.

As above by shrinking $\epsilon$ if necessary, we can assume the image  $F(D_{\epsilon})$
is contained in some coordinate chart $(U, w)$ of $X_0$.
Now consider the inverse of $F|_{N_{\epsilon}}$, which we denote by $H: \widetilde{N} \rightarrow
 N_{\epsilon} \subset \mathbb{C}^n.$ By the Hartogs extension, $H$ extends to a holomorphic function in
  $\widetilde{D},$ which we still denote by $H.$  By the uniqueness of holomorphic maps and by
  the maximum principle, $H(\widetilde{D})$ is contained in $D_{\epsilon}.$ Furthermore,
   we have $\hat{F} \circ H$ equals the identity map on $\Ol{\widetilde{D}}$. (In particular, $H$ is embedding). On the other hand, we also have $H \circ F$ equals identity on $\Ol{D_{\epsilon}}-\{ p\}$ and by continuity, it also equals identity on $\Ol{D_{\epsilon}}.$ This implies $H$ is a biholomorphic map from $\widetilde{D}$ to $D_{\epsilon}$.
 Hence we prove $\mathrm{sing}(\Omega)$ is empty.

Once we know $\Omega$ has no singularities and thus it is a complete K\"ahler manifold equipped with the Bergman metric,  a standard argument will show that $\Omega$ must be biholomorphic to the ball. Indeed, we now have a holomorphic isometric immersion $F$ from $\Omega$ to $X_0$. On the other hand, we
consider the local inverse $(g, V)$ of $(f, U)$, where $V=f(U).$ It is a local isometric embedding from $f(U)$ to $\Omega.$ Since now $(\Omega, \omega^B)$ is complete,  $(g, f(U))$ extends holomorphically along any path in $X_0$ (See Proposition 11.3, 11.4 in \cite{He}). As $X_0$ is simply connected,  we thus obtain a holomorphic map $G$ from $X_0$ to $\Omega$. Note $f \circ g$ equals the identity map on $V$ and $g \circ f$ equals the identity map on $U$. By the uniqueness of holomorphic functions, $F \circ G$ and $G \circ F$ equal to the identity map on $X_0$ and $\Omega$, respectively. We thus conclude $\Omega$ is biholomorphic to $X_0.$ Finally, it follows from the assumption on $\Omega$ and $\omega^B$ that $X_0$ can only be the complex unit ball. Furthermore, as the holomorphic sectional curvature of the Bergman metric is a biholomorphic invariance, we have $\lambda=-\frac{2}{n+1}.$ \endpf

\end{document}